\documentclass{amsart}
 \usepackage{graphicx}

 \usepackage{graphicx}
 \usepackage{epstopdf}
     \usepackage{amssymb}
  \usepackage{color} 
 \usepackage{amsthm,graphpap}

    \newcommand\nilp{{\mathrm{nilp}}}

\long\def\symbolfootnote[#1]#2{\begingroup%
\def\thefootnote{\fnsymbol{footnote}}\footnote[#1]{#2}\endgroup}

\newcommand\CD{\check{D}}

\newcommand\fbul{F^\bullet}
\newcommand\fbulp{{F'}^{\bullet}}
\newcommand\fbuly{\fbul_{\infty}}
\newcommand\fbulim{\fbul_{\lim}}
\newcommand\tfbul{\wt F^\bullet}
\newcommand\bbul{{\bullet,\bullet}}
\newcommand\Obul{\Omega^\bullet}

\newcommand\bsm{ \begin{smallmatrix}}
\newcommand\esm{\end{smallmatrix} }
\newcommand\bbm{\left[\begin{smallmatrix}}
\newcommand\ebm{\end{smallmatrix}\right]}
\newcommand\bcs{\begin{cases}}
\newcommand\ecs{\end{cases}}
\newcommand\wbul{{W_\bullet}}
 
\newcommand{\lrc}[1]{\left\{ #1\right\}}

\newcommand{\lrp}[1]{\left(#1\right)}
 \newcommand\bp[1]{\bigl(#1\bigr)}
  \newcommand\Bp[1]{\Bigl(#1\Bigr)}

 \newcommand\wt[1]{\widetilde{#1}}
  \newcommand\wh[1]{\widehat{#1}}

 \newcommand\sidecong{\rotatebox{90}{$\cong$}}
 
  \newcommand\sideeq{\rotatebox{90}{$=$}}

\newcommand{\C}{{\mathbb{C}}}

\newcommand\bH{\mathbb{H}}

\renewcommand{\P}{\mathbb{P}}
\newcommand{\Q}{{\mathbb{Q}}}

\newcommand{\R}{\mathbb{R}}

\newcommand{\Z}{\mathbb{Z}}

\newcommand{\cC}{{\mathscr{C}}}

\newcommand{\cF}{{\mathscr{F}}}
\newcommand{\cG}{\mathfrak{g}} 
\newcommand\cGC{{\cG_\C}}
\newcommand\cGR{{\cG_\R}}

\newcommand{\cM}{{\mathscr{M}}}
\newcommand{\cN}{{\mathscr{N}}}
\newcommand{\cO}{{\mathscr{O}}}
\newcommand{\co}{{\scriptstyle\mathscr{O}}}

\newcommand{\cS}{{\mathscr{S}}}

\newcommand{\cX}{{\mathscr{X}}}
\newcommand{\cY}{{\mathscr{Y}}}

\newcommand\bpm{\begin{pmatrix}}
\newcommand\epm{\end{pmatrix}}

\newcommand\mtd{Mumford-Tate domain}

\newcommand\phs{polarized Hodge structure}
\newcommand\mhs{mixed Hodge structure}

\newcommand\lmhs{limiting mixed Hodge structure}

\renewcommand\mod{\mathop{\rm mod}\nolimits}

\newcommand\rspan{\mathop{\rm span}\nolimits}

\newcommand\rsl{\mathop{\rm sl}\nolimits}

 \newcommand{\Gr}{\mathop{\rm Gr}\nolimits}

\newcommand\End{\mathop{\rm End}\nolimits}

\newcommand\rim{\mathop{\rm Im}\nolimits}

\newcommand\sing{{\mathop{\rm sing}\nolimits}}

\newcommand\prim{{\rm prim}}
\newcommand\Res{\mathop{\rm Res}\nolimits}

\newcommand{\Ext}{\mathop{\rm Ext}\nolimits}

 \newcommand{\Hom}{\mathop{\rm Hom}\nolimits}

\newcommand{\AJ}{{\rm AJ}}

 \renewcommand{\part}{\partial}
\newcommand{\la}{{\lambda}}
\newcommand\ga{{\gamma}}
\newcommand\Ga{{\Gamma}}

\newcommand{\Om}{{\Omega}}
\newcommand{\om}{{\omega}}

\newcommand\vp{\varphi}
\newcommand\sig{\sigma}

\newcommand{\chs}{\hskip-10pt}

\newcommand\bsl{\backslash}
\newcommand\eps{\epsilon}

 \renewcommand{\AA}{{\mathfrak{A}}}

\newcommand{\lab}{\label}
\newcommand\pref[1]{{\rm (\ref{#1})}}

\newcommand{\hensp}[1]{\enspace\hbox{#1}\enspace}

\newcommand{\opplus}{\mathop{\oplus}\limits}

  \renewcommand\thesection{\arabic{section}}

\newcounter{demo}[equation]

 \newenvironment{demo}{
 \addtocounter{equation}{1}\refstepcounter{demo} \setcounter{demo}{\value{equation} }
{\smallbreak\noindent (\thesection.\arabic{demo})\enspace   }} 
 
 \newtheoremstyle{mytheo}
 {3pt}
  {3pt}
  {\itshape}
  {}
  {\scshape}
  {:}
  {.5em}
  {}

\theoremstyle{mytheo}

\newtheorem{thm}{Theorem}

\newtheorem{thmp}{Theorem}

\newtheorem*{propn}{Proposition}

\newtheoremstyle{note}
  {3pt}
  {3pt}
  {}
  {}
  {\bfseries}
  {:}
  {.5em}
  {}
\theoremstyle{note}

\newtheorem*{exam}{Example}

\theoremstyle{remark}

\newcommand\simto{\xrightarrow{\sim}}
\newcommand\xri[1]{\xrightarrow{#1}}

\newcommand\bx{\mathbf{x}}
\newcommand\bt{\mathbf{t}}

\newcommand\ol[1]{\overline{#1}}

 \usepackage[mathscr]{euscript}
\input{xy}
\xyoption{all}
  \renewcommand\thesection{\arabic{section}}

  \usepackage{epic}
  \renewcommand\thesection{\Roman{section}}
 \let\pr\pref

\newcommand\Def{\mathop{\rm Def}\nolimits}
\newcommand\Gy{\mathop{\rm Gy}\nolimits}
\newcommand\Rest{\mathop{\rm Rest}\nolimits}
\newcommand\EExt{\mathop{\mathbb{E}\mathrm{xt}}\nolimits}
 
 \newcommand\Spec{\mathop{\rm Spec}\nolimits}
 \newcommand\plmhs{polarized \lmhs}
 \newcommand\lpm{limit period mapping}
 \newcommand\cOX{{\cO_X}}
 \newcommand\cxxi{{\cX_\xi}}
  \newcommand\cOD{{\cO_D}}
  \newcommand\cOXt{{\cO_{\wt X}}}
 \newcommand\Omox{\Om^1_X,\cOX}
  
 \newcommand\es{{\rm es}}
 
\newcommand\fst{{\ensuremath$1^{\rm st}$} order}
\newcommand\vspth{\vspace*{-3pt}}
\newcommand\bxi{{\boldsymbol{\xi}}}
\newcommand\bla{{\boldsymbol{\la}}}
 \newcommand\cxxib{{\cX_\bxi}}
\newcommand\lms{limiting mixed structure}
\newcommand\beps{{\boldsymbol{\eps}}}
\newcommand\abul{A^\bullet}
\begin{document}

\title[Deformation theory]{Deformation theory and limiting mixed Hodge structures}
 \author[M. Green and P. Griffiths]{Mark Green and Phillip Griffiths}
 \address{Department of Mathematics, University of California at\hfill\break\indent Los Angeles, Los Angeles, CA 90095\hfill\break\indent
 {\it E-mail address}\/: {\rm mlg@ipam.ucla.edu}\hfill\break\indent  Institute for Advanced Study, Einstein Drive, Princeton, NJ 08540\hfill\break\indent
{\it E-mail address}\/: {\rm pg@ias.edu}}
      \maketitle

     \section*{Outline}
     \begin{itemize}
     \item[I.] Introduction
     \item[II.] Deformation theory
     \item[III.] Proof of Theorems I and I$'$
     \item[IV.] Proof of Theorem II
     \item[V.] Proof of Theorem III 
     \item[VI.] The hierarchy of mixed Hodge structures
     \item[] References
     \end{itemize}
     \vspace*{-8pt}
     \section{Introduction}
     
This paper was  motivated by the following question: Recall that for a smooth projective variety $X$ whose \phs\  on $H^n(X,\Q)_\prim$ leads to  a period point $\Phi(X)\in D$, the period domain for \phs s of a fixed type, the differential
     \vspth\[
     \Phi_\ast: T_X\Def(X)\to T_{\Phi(X)}D    \vspth\]
     from the tangent space to the Kuranishi space $\Def(X)$ to the tangent space  $T_{\Phi(X)}D$ is expressed cohomologically by the inclusion $T_X \Def(X)\allowbreak\hookrightarrow \EExt^1_{\cO_X}\bp{\Om^1_X,\cO_X}\allowbreak\cong H^1(T_X)$,
  and the resulting natural maps on the associated graded to  $\fbul \bH^m(\Obul_X)$  
    induced from
       \vspth  \begin{equation}\lab{1.1}
     \EExt^1_{\cO_X} \bp{\Om^1_X,\cO_X}\otimes \bH^m(\Obul_X)\to \bH^m(\Obul_X). \footnotemark    \vspth\end{equation}
     \footnotetext{In more classical notation   $\Theta_X=T_X$, this is the map
     \[
     H^1 \lrp{\Theta_X}\otimes H^{m-p}\lrp{\Om^p_X}\to H^{m-p+1} \lrp{\Om^{p-1}_X}\]
     induced from the contraction $\Theta_X\otimes \Om^p_X\to \Om^{p-1}_X$.}
Our question was how  to extend this formalism to the case where $X$ is singular, having singularities of the type that arise by semi-stable reduction in a family of projective varieties $\cX\xri{f} S$ whose general member is smooth.  This question was   studied by Friedman \cite{Fr1} when $\dim S=1$.  We were interested in amplifying and extending his results, emphasizing the development of a formalism that lends itself to the computation of examples.  In the course of trying to carry  this out we have found that there is an interesting story surrounding the relationships among the various \mhs s  associated to $X$ and   its \fst\ neighborhood in $\cX$, and one of the purposes of this paper has turned out to be to amplify and clarify these relationships in the context of deformation theory.
  Here for us  the works  \cite{Fr2}, \cite{St1} and \cite{St2} have been very important  when $\dim S=1$, as has the extension  of \cite{St1} to a general $S$ by Fujisawa \cite{Fu1}, \cite{Fu2}.  In fact, this  is a partly expository paper, drawn from, reinterpreting and building on the works \cite{Fr1}, \cite{Fr2}, \cite{St1}, \cite{Zu},  \cite{PS}, \cite{St2}, \cite{Fu1}, \cite{Fu2}, \cite{CKS1}, \cite{KP2},
         \cite{GGR} and others. 
     
     To address the question stated above one is led to focus on the singular variety $X$ and its first order deformations.  For the case when $X$ is a normal crossing variety this is done in \cite{Fr2}.  Here motivated by the semi-stable reduction theorem in \cite{AK} we shall assume more generally that 
     \[
     \hbox{\em $X$ is locally a product of normal crossing varieties.}\]
     This means locally in $\C^n$ with coordinates $x_1,\dots,x_{n}$ we have a sequence $1\leqq i_1<i_2<\dots<i_k\leqq n$ with index blocks $I_1=\{1,\dots,i_1\}$, $I_2=\{i_1+1,\dots,i_2\},\dots$, and then  $X$ is locally given by
     \begin{equation}\lab{new1.2}
    x_{I_1}=0,\dots,x_{I_k}=0\end{equation}
     where $x_{I_1}=x_1\cdots x_{i_1}$, $x_{I_2} = x_{i_1+1}\cdots x_{i_2},\ldots$\,.  The usual locally normal crossing variety is the case $k=1$.\footnote{There is an important distinction between the case when $X$ is locally a normal crossing variety and when it is globally such.  By a combination of  blowings up and base changes  the former may be reduced to the later,  and for both theoretical and notational purposes this is generally done.  For computational purposes the former is frequently more convenient; e.g., for irreducible nodal curves.  In this paper we shall restrict to the global normal crossing case and its analogue when $X$ is locally a product of normal crossing varieties.  However, we expect that the discussion given below  will remain valid in the more general case, and some of our examples are carried out in the local normal crossing case. The formalism in \cite{De} and \cite{St2} allows one to handle the general theory when $X$ is locally a  normal crossing variety.}  The deformation theory of such varieties is well understood \cite{Pa}, and for simplicity of exposition in this paper we shall abuse notation and set
     \[
     T_X \Def(X)=\EExt^1_{\cO_X} \bp{\Om^1_X,\cO_X}.\]
     The abuse of notation is because here the right-hand side is the  space of deformations of $X$ over $\Delta_\eps=: \Spec \C [\eps]$, $\eps^2=0$, so that it is only the Zariski tangent space to the Kuranishi space $\Def(X)$.  
  In general there may  be obstructions to lifting deformations defined over the Zariski tangent space, but this issue will play no role in what follows.\footnote{One may make the blanket assumption that all \fst\ deformations are unobstructed, and then at the end note that this assumption has never   been used.  }   In fact, one of the main points is that the theory of \lmhs s for $1$-parameter families depends only on the $1^{\rm st}$ order neighborhood of the singular variety, a point that is     implicit in \cite{Fr2} and explicit in a somewhat different form  in \cite{St2}.

     We shall   make the crucial assumption that \emph{there exists a $\xi\in\EExt^1_{\cOX} \bp{\Om^1_X,\cOX}$ such that for every $x\in X$ the  image $\xi_x$ of $\xi$ in the natural map}
    \begin{equation}\lab{1.2}
    \EExt^1_{\cO_X}\bp{\Om^1_X,\cO_X} \to \Ext^1_{\cO_X} \bp{\Om^1_X,\cO_{X}}_x\end{equation}
\emph{smooths to \fst\     the singularity at $x$.}  Equivalently, for every $x\in X$  the global deformations of $X$ over $\Delta_\eps$ map to  smoothing deformations 
      of the germ $X_x$ of $X$ at $x$.  The smoothing deformations of  \pr{new1.2} are given by $x_{I_j}=t_j$ and they have tangents $\sum_{i=1}\la_i \part_{t_i}$ where all $\la_i\ne 0$. We denote by
    \[
    T^0_X \Def(X)\subset T_X\Def(X)\]
    the open set of all $\xi\in T_X\Def(X)$ whose localizations are smoothing deformations of $X_x$ for every $x\in X$.
    
    We define the pair $(X,\xi)$ to be \emph{projective} in case there is a very ample line bundle $L\to X$ such that $L$ extends to $\cX_\xi$. This can be expressed cohomologically in a standard way, and we shall assume it to always be the case.   
    
    A \emph{\lmhs}\ $(V,\wbul,\fbul)$ is given by a $\Q$-vector space $V$, a weight filtration $\wbul$ and Hodge filtration $\fbul$ that define a \mhs, and where there exists a nilpotent $N\in \End(V)$ such that (i) $\wbul = \wbul(N)$ is the monodromy weight filtration, and (ii) for the integer $m$  around which the monodromy weight filtration is centered, the $N^k:\Gr^{\wbul}_{m+k}V \simto \Gr^{\wbul}_{m-k} V$ are isomorphisms for $0\leqq k\leqq m$.  The \lmhs\ is \emph{polarizable} if there exists a $Q:V_\otimes V\to \Q$ with $Q(u,v)=(-1)^m Q(v,u)$, and an $N\in\End_Q(V)$ as above such that the \emph{primitive spaces} $\ker N^{k+1}$ above  are polarized Hodge structures via $Q_k(u,v)=\pm Q(N^k u,v)$ (cf.\ \cite{CKS1}).  Two \lmhs s $(V,\wbul,\fbul)$ and $(V,W'_\bullet,\fbulp)$ are \emph{equivalent} if $W'_\bullet = \wbul$, and if $\fbulp=\exp(zN)\fbul$ for some $z\in \C$.  We will denote by $[V,\wbul,\fbul]$ an equivalence class of \lmhs s.
    
We shall   use the term \emph{standard family} to mean that $\cX_\Delta\to \Delta$ is a projective mapping where $\cX_\Delta$ is smooth, the fibres $X_t=\pi^{-1}(t)$ are smooth for $t\neq 0$, and $X_0=X$ is a reduced   normal crossing variety.
      
    \begin{thm} \lab{thm1}
    Canonically associated to each $\xi\in T^0_X\Def(X)$ is a \lmhs\ $(V_\xi,\wbul,\fbul_\xi)$.  In case
    $X$ is a   normal crossing variety and 
     $\xi$ is tangent to an arc $\Delta\subset \Def(X)$ giving a standard  family $\cX_\Delta\xri{\pi}\Delta$ with $\pi^{-1}(0)=X$, this \lmhs\ is the one associated to the family and $\xi\in T_0(\Delta)$.
    \end{thm}

    This result is largely an amalgam and slight extension of those in   \cite{Fr2} and \cite{St2}.
    A key point is to note that the data $(X,\xi)$ gives  a standard family    $\cX_\xi\to \Delta_\eps$, together with  an extension
    \begin{equation}\lab{1.3}
    0\to\cOX\to \Om^1_{\cX_\xi}\otimes \cOX\to \Om^1_X\to 0\end{equation}
of $\cOX$-modules.\footnote{We are here extending the notion of a standard family to include the smooth non-reduced scheme $\cxxi$ with structure sheaf $\cO_{\cxxi}$ locally isomorphic to $\cOX[\eps]$.  We   will also say that fibres over $\Delta^\ast_e$ are smooth.}      A second key point is to show that, as explained in section III below, \pr{1.3}     gives   an exact sequence
  \begin{equation} \lab{1.5}
  0\to \cOX\to \Om^1_{\cX_\xi}(\log X)\otimes \cOX\to \Om^1_{\cX_\xi/\Delta_\eps}(\log X)\otimes \cOX\to 0,\end{equation}
  where in the case of a   standard family $\cX_\Delta\to \Delta$  restricting to $\cxxi\to \Delta_\eps$
  \[
  \Om^1_{\cX_\xi/\Delta_\eps}(\log X)\otimes \cOX=\Om^1_{\cX_{\Delta/\Delta}}(\log X)\otimes \cOX.\]
  The vector space in the \lmhs\ is   
  \[
  V_\xi = \bH^m\lrp{\Obul_{\cX_\xi/\Delta_\eps} (\log X)\otimes \cOX}\]
  and $\fbul_\xi$ is induced from the ``b\^ete'' filtration   on $\Obul_{\cX_\xi/\Delta_\eps} (\log X)\otimes \cOX$.  The monodromy logarithm is induced from the connecting homomorphisms arising from \pr{1.5}. The $\Q$-structure and properties of the  monodromy logarithm and resulting monodromy weight filtration  are more subtle to define and treat  (cf.\ \cite{St1}, \cite{Zu} and chapter 11 in \cite{PS}).

  We note that the usual ambiguity in either the Hodge filtration or the $\Q$-structure in the \lmhs\ associated to $\cX\to \Delta$, ambiguity   that depends on a choice of parameter $t$, is removed by considering the data $(X,\xi)$.
  
  A subtle point, one that will be further explained below, is this: For $X$ a smoothable normal crossing variety the singular locus $D$ will have connected components $D_a$.  Then we will see that  $\Ext^1_\cOX(\Omox)\cong \oplus \cO_{D_a}$, and in the basic exact sequence \pr{2.1} a global \fst\ deformation $\xi\in\EExt^1_\cOX(\Omox)$ will induce $\xi_{D_a}\in H^0(\cO_{D_a})$.  The condition that $\xi$ be to \fst\ smoothing along $D_a$  is that $\xi_{D_a}\ne 0$.  \emph{Then the equivalence class of the   \lmhs\ in Theorem \ref{thm1} depends only on the $\xi_{D_a}$ and not on the global $\xi$ that maps to the $\xi_{D_a}$'s.}  In fact, given a collection  of non-zero $\xi_{D_a}$'s, we may construct a \lmhs\ provided that there is a global smoothing $\xi$; the particular $\xi$ does not matter.\footnote{As we hope to discuss further in a work in progress, this is related to the theorem of Cattani-Kaplan \cite{CK} that the weight filtration $\wbul(N)$ is independent of $N$ in the interior of a monodromy cone, and the result in \cite{CKS1} that the equivalence class of the  \lmhs\ is independent of the direction of approach from the interior of the cone.}
  
  To handle several variable families we shall consider a vector
  \[
  \bxi\in \EExt^1_\cOX\bp{\Om^1_X,\cOX^\ell} = \opplus^\ell T_X\Def(X)\]
  with the property that for $\la=(\la_1,\dots,\la_\ell)$ and 
  \[
  \bxi_\la=\la_1\xi_1+\cdots+\la_\ell\xi_\ell,\qquad \la_i\ne 0\]
  we have
  \begin{equation}
  \lab{1.6}
  \bxi_\la\in T^0_X\Def(X).\end{equation}
  For $\Delta_\beps = \Delta_{\eps_1}\times \cdots \times \Delta_{\eps_\ell}$, it will be seen that  we then have a family
  \[
  \cX_{\bxi_\la}\to \Delta_\beps\]
  with   smooth fibres over $\Delta^\ast_\beps = \Delta^\ast_{\eps_1}\times \dots\times \Delta^\ast_{\eps_\ell}$.  We think of this as a space of \fst\ deformations that deform $X$ to a ``less singular" variety along the axes but which smooth  $X$ when we deform into the interior.\footnote{In the paper \cite{KN} the definition of a normal crossing variety with logarithmic structure is introduced.  The presence of a logarithmic structure is equivalent to $d$-semi-stability in the sense of \cite{Fr2} 
    (cf.\ \pr{2.6} below).  A deformation theory for normal crossing varieties with logarithmic structure is then introduced.  In the context of  this paper this theory amounts to deformations of $X$ that independently smooth the connected components of the singular locus $D$ of $X$, modulo equisingular deformations.  The  log-geometry formalism nicely lends itself to computation of examples for Calabi-Yau varieties.}   
     There is then a  several variables  analogue of Theorem \ref{thm1} where the terms in the statement will be explained in the text.
  
  \begin{thmp} \lab{thm1p}
  Associated to $\bxi\in\EExt^1_\cOX\bp{\Om^1_X,\cOX^\ell}$ satisfying the condition that  \pr{1.6} holds, there is a several variable \lmhs\ $(V_\bxi,\wbul,\fbul_\beps)$ in the sense of \cite{CKS1}.  In case $X$ is a   normal crossing variety and $\bxi_\la$ is tangent to an arc $\Delta_\la\subset \Def(X)$, this \lmhs\ is the one associated to the standard family $\cX_{\Delta_\la}\xri{\pi}\Delta_\la$ with $\pi^{-1}(0)=X$.\end{thmp}
  
  As will be discussed below, for the last statement in the theorem the general case when $X$ is locally of the form \pr{1.2} seems to be open (cf.\ \cite{Fu1}, \cite{Fu2}), and we will discuss a geometric reason for this.
  
  Detailed proofs of Theorems I and \ref{thm1p}, especially for the latter, will not be given below.  The argument for Theorem \ref{thm1} consists largely of proof analysis of those in the references \cite{Fr2}, \cite{St1} and \cite{St2} and will be addressed more fully in a work in progress.  For Theorem \ref{thm1p}, the construction of a \mhs\ follows largely from our construction given below and \cite{Fu1}, \cite{Fu2}.  The construction of a \emph{limiting} \mhs\ requires more work  and will be taken up in the work in progress. We will however try to point out some of the key points in both of these arguments.
  
  For the analogue of \pr{1.1} we have
  \begin{thm}\lab{thm2}
  The class $\xi\in\EExt^1_{\cOX}\bp{\Om^1_X,\cOX}$ in \pr{1.3} defines a natural class $\xi^{(1)} \in \EExt^1_{\cOX}\bp{\Om^1_{\cX_\xi/\Delta_\eps}(\log X)\otimes \cOX,\cOX}$ corresponding to \pr{1.5}, and the $1^{\rm st}$ order variation of the \lmhs\ is expressed as the natural mapping
\[\displaylines{\qquad 
\EExt^1_\cOX\lrp{\Om^1_{\cxxi/\Delta_\eps} (\log X)\otimes \cOX,\cOX}\hfill\cr\hfill \to \End_{\rm LMHS} \bH^m
\lrp{\Obul_{\cxxi/\Delta_\eps}(\log X)\otimes \cOX}. \footnotemark\qquad}\]
  \end{thm}
  \footnotetext{The notation $\xi^{(1)}$ has been used because the construction of the sequence \pr{1.5} from \pr{1.3} resembles that of the construction of the first prolongation in the theory of exterior differential systems.  The group $\EExt^1_\cOX \bp{\Om^1_{\cX_\xi/\Delta_\eps} (\log X)\otimes \cOX,\cOX}$, which may be defined if there exists a logarithmic structure on $X$, appears naturally in the deformation theory of smooth logarithmic varieties (cf.\ \cite{Abetal}).
  
  Referring to footnote \ref{17} below, in the setting of log-analytic geometry the important monograph \cite{KU} contains a treatment of the differential of the period map at infinity for standard families $\cX_\Delta\to\Delta$ (cf.\ Theorem 4.4.8).  In case the $\xi$ in Theorem II arises as the tangent vectors at the origin we believe that those results should be equivalent.}
  
  \noindent Here, $\bH^m \lrp{\Obul_{\cxxi/\Delta_\eps}(\log X)\otimes \cOX}= V_\xi$ is the vector space underlying the \lmhs\ in Theorem \ref{thm1}, and $\End_{\rm LMHS}$ means the endomorphisms of $V_\xi$ that preserve the structure as a \lmhs\ as explained below.
  
  Again the terms in the statement will be explained in the text.  An informal way to think about this result is this: Denoting by $\Def(X ,\xi)$ the deformations of the pair $(X,\xi)$, we have a natural \emph{extended period mapping}
  \[
  \Def(X,\xi)\to \CD\]
  that assigns to $\xi\in T^0_X\Def(X)$ the well-defined point $\fbul_\xi\in\CD$, the dual space to the period domain $D$   consisting of filtrations  of $V_\xi$ that satisfy  only the $1^{\rm st}$ Hodge-Riemann bilinear relation.  Then in   the map in Theorem \ref{thm2} might be thought to be the differential
  \begin{equation} \lab{new1.7}
  T_{(X,\xi)} \Def(X,\xi)\to T_{F^\bullet_\xi}\CD\end{equation}
  of the extended period mapping.  This is \emph{not} the case, as will be   made   precise in Section IV below. 
   The issue is more subtle in that $\xi$ gives not only a well-defined \lmhs,  not just an equivalence class of such, but also defines a \fst\ variation of that \lmhs.  This is the information in $\xi^{(1)}$.   At first glance one might think that since it takes the tangent vector $\xi$ to define $F^\bullet_\xi$, the information in $\xi^{(1)}$ which gives the variation of  the entire \lmhs\ would be of $2^{\rm nd}$ order.  But this is not correct, and it was in trying to understand this loint that we were led to most of the other topics in this paper.

   We will however  see by example that $\xi^{(1)}$ contains strictly more information than the differential at the origin of the Kato-Usui map \cite{KU}
  \[
  \Delta\to \Ga_T\bsl D_N.\]
Here $D_N=D\cup B(N)$ is the period domain $D$ with the boundary component $B(N)$  attached to $D$, where $B(N)$ consists of all equivalence classes of \lmhs s with monodromy logarithm $N$   and where   $\Ga_T = \{T^\Z\}$ with $T=\exp N$ is the local monodromy group  
  (cf.\ Section IV below for an explanation of the notations and terms used).  It is in this sense that Theorem \ref{thm2} provides an answer to our original question.  The term ``expressed" means that in examples $\EExt^1_\cOX\bp{\Om^1_{\cxxi/\Delta_\eps}(\log X)\otimes \cOX,\cOX}$ will have algebro-geometric meaning    and the pairing is a   cup-product.  We will see by example that the additional information is non-trivial and somewhat subtle.\footnote{In very classical terms one may write the period matrix $\Om(t)$ in block form where the blocks $\Om_i(t)$ are polynomials in  $\log t$ with holomorphic coefficients and where the remaining blocks $\Om_\alpha(t)$ are holomorphic at $t=0$.  The differential of the map to $\Ga_T\bsl D_N$ records the derivatives $\Om'_\alpha(0)$ of the holomorphic terms, while \pr{new1.7} has the effect of regularizing the logarithmically divergent integrals that give the $\Om_i(t)$ and then taking the linear part $\Om'_i(0)$ at $t=0$ of that regularization.  The $\Om'_\alpha(0)$ and $\Om'_i(0)$ record the variation in the full extension data in the \lmhs.}
  
  As will be explained  in Section  V below, associated  to a \plmhs\ is a \emph{reduced \lpm}\ and distinguished point
  \[
  \fbul_\infty \in\part D \]
  where $D= G_\R/H$ is a period domain (\cite{KP1}, \cite{KP2}, \cite{GGK},  \cite{GG} and \cite{GGR}).
  The boundary $\part D$ is stratified into finitely many $G_\R$-orbits and their geometry is a much studied and very interesting topic (\cite{GGK}, \cite{FHW}).
  
  On the other hand, the vector space $T_X\Def(X)$ is stratified by open sets $T^0_X\Def(X)_I$ contained in linear subspaces $T_X\Def(X)_I\subset T_X\Def(X)$.  In the text we will explain this in case $X$ is a normal crossing divisor, which is the only case    for which  thus far we have a result.  Then the strata correspond to subsets of the set of connected components of the singular locus $X_\sing$ of $X$.  The subspace $T^0_X\Def(X)$ is the open stratum of smoothing deformations; the other strata correspond to the components that are smoothed when $X$ deforms in the directions of that strata.  The opposite extreme to $T^0_X\Def(X)$ is the linear subspace $T_X\Def^{\es}(X)\cong H^1\lrp{\Ext^0_{\cOX}\bp{\Om^1_X,\cOX}}$ of equisingular deformations.
  It seems reasonable to expect,  but     we are not aware of  a proof  in the literature, that $\xi\in T^0_X\Def(X)_I$ corresponds to the limit in a variation of \mhs s over the punctured disc (\cite{St-Zu}).
  
  Leaving this important issue aside, we return to the deformation theory and \lmhs s in the several parameter case.
      In the study of \lmhs s over higher dimensional base spaces (\cite{CKS1}) there are a number of cone  structures that enter:
  \begin{enumerate}
  \item the stratification of abelian subspaces $\AA\subset \cG^{\nilp}$ induced by the $G$-orbit structure on $\cG^\nilp$ (\cite{Ro} and references cited therein);
  \item the stratification of nilpotent cones as in \cite{CKS1} and \cite{KU} (cf.\ \cite{AMRT} for the classical weight one case);
  \item the stratification of $\part D$ by $G_\R$-orbits and its relation to reduced \lpm s  \cite{KP1}, \cite{KP2}, \cite{GG},   \cite{GGR} and \cite{Ro} and work in progress by Kerr, Pearlstein and Robles; and
  \item the stratification of $T_X\Def(X)$, as explained below for $X$ a normal crossing divisor, and which  we feel can reasonably be   expected  to extend to the case where $X$ is locally a product of normal crossing divisors.\end{enumerate}
  The basic known result, due to Robles \cite{Ro}, is that \emph{the interiors  of the strata in (ii) map to strata in (i)}, and as a consequence to strata in (iii).\footnote{For the interior of the full nilpotent cone this result follows from \cite{CKS1}.}  Her argument makes full use of the deep properties of several variable nilpotent orbits \cite{CKS1} and of the classification of $G_\R$-orbits in $\cG^{\nilp}_\R$ (cf.\ the references in \cite{Ro}).  An algebro-geometric version of Robles' result might be  that at the tangent space level strata in (iv) map to strata in (iii).  The theorem to be described now  is a partial result in this direction.
  
  In the setting of the Cattani-Kaplan-Schmid theory there are defined nilpotent cones
  \[
  \sigma=\rspan_{\Q> 0} \{N_1,\dots,N_\ell\}\]
  where the $N_i\in\cG^\nilp$ are linearly independent commuting nilpotent transformations and several variable nilpotent orbits $(\fbul,\sig)$.
  Here, $\fbul\in\CD$ and the conditions
  \begin{itemize}
  \item $\exp(z_1N_1+\cdots +z_\ell N_\ell)\cdot \fbul\in D$ for all $\rim z_i\gg 0$;
  \item $[N_i,F^p]\subset F^{p-1}$
  \end{itemize}
  are satisfied.  We denote by $\wt B(\sig)\subset \CD$ the set of several variable nilpotent orbits, and by $B(\sig)$ the equivalence classes of those orbits under reparametrization $z_i\to z_i+\la_i$.  In \cite{KP1}, \cite{GGK} and \cite{GG}   there are defined reduced \lpm s  for 1-dimensional cones, and the construction can be extended \cite{KP2} to the general  case      to give the reduced limit period map
  \begin{equation}\lab{new1.8}
  \Phi_\infty: B(\sig)\to \part D.\end{equation}
  
  \begin{thm}\lab{thm3}
  Let $X$ be a normal crossing variety for which there exists a  $\xi\in T_X\Def(X)$ that is nowhere vanishing along each component of $X_\sing$.  Then there exists a nilpotent cone
  \[
  \sig_X\subset T_X\Def(X)/T^{\es}_X\Def(X),\]
  and a several variable \lmhs\ in the sense of \cite{CKS1} with the property that under the reduced \lpm\ \pr{new1.8} $\wt B(\sig_X)$ maps to a $G_\R$-orbit in $\part D$.\end{thm}
  This theorem follows from the construction of $\sig_X$ and the result of Robles mentioned above. As mentioned before, it is of interest to see if the construction of $\sig_X$ and the result can be extended to the faces of the cone $\sig_X$.
  
  For our next result we note that given a standard family $\cX\to \Delta$   there are the following four types of mixed Hodge structures that may be defined:
  \begin{enumerate}
  \item the part of the mixed Hodge structure on $H^\ast(X)$ that comes from the \lmhs;\footnote{This is $\ker N$.}
  \item that part of the \lmhs\ that may be defined in terms of $X$ alone;
  \item the limiting mixed Hodge, modulo reparametrizations $\fbulim\sim  \exp(z N)\cdot \fbulim$ resulting from a change in parameter in the disc, associated to $\cX\to \Delta$;\footnote{This is by definition the same as an equivalence class of \lmhs s.} and
  \item the \lmhs\ associated to $(X,\xi)$, where $\xi\in T^0_X\Def(X)$ is the first order variation of $X$ in $\cX$.\end{enumerate}
  
  \begin{thm} \label{th4}
  There are strict implications
 \[ {\rm (iv)}\implies {\rm (iii)}\implies {\rm (ii)}\implies {\rm (i)}.\]
   \end{thm}
   The term ``strict implication'' means  that there is successively more information in (i), (ii), (iii), (iv); the precise meaning of this will be explained in the proof.
   
    We will see that given an abstract $X$  that is locally a normal crossing divisor, the condition that we may construct the data given in (ii) is that there exists a $\xi\in T^0_X \Def(X)$ that is smoothing to \fst;  the actual data will not depend on the particular $\xi$ but rather will depend on the $\xi_{D_a}$'s as discussed above.  A \lmhs\ will decompose into $N$-strings under the action of the monodromy logarithm $N$.  This decomposition may be pictured as
 \begin{equation}    \lab{1.9} 
  \begin{matrix}  
  H^0(-m) \longrightarrow \hspace*{4pt}  \quad \qquad\cdots\hspace*{-4pt}  \qquad \longrightarrow H^{0}(-1) \longrightarrow H^0 \\
    H^{1} (- m+1)   \longrightarrow  \cdots\longrightarrow  H^{1}(-1)\longrightarrow  H^{1}\\
   \vdots\quad\ \ \\
    H^m\quad \  \end{matrix}\end{equation}
  where $H^k$ is a pure Hodge structure of weight $k$.\footnote{If one thinks of $N$ as being completed to an $\rsl_2$-triple, then the $N$-strings are composed from the irreducible pieces in the decomposition of the $\rsl_2$-module.  The $H^k$'s on the right end may themselves be Tate twists of lower weight Hodge structures.}   We may think of \pr{1.9} as giving the primitive decomposition in the associated graded to a \lmhs, together with the iterated action of $N$ on the primitive spaces.  Then our result pertaining to (ii) is

  \begin{thm}\lab{thm6}
  The terms $H^{m-j}(-i)$, $0\leqq i\leqq m-j$ in \pr{1.9},  together with the $N$-maps between them,  may be constructed from $X$ alone.\end{thm}

  We will also see for $[\xi]\in\P T^0_X \Def(X)$  with localizations $\xi_{D_a}$ along the components of $D_a$ of $X_\sing$,  we will have 
  \[
  {\rm (iii)}\longleftrightarrow (X,[\xi_{D_a}]\hbox{'s}) ,\]
  and where  the brackets refer to the corresponding point in the designated projective space  and the symbol ``$\longleftrightarrow$" means that the data on each side are equivalent.
  
  We hope that this result will clarify exactly what input is needed to be able to define the limiting mixed Hodge structures, or the parts thereof, that are associated to a degeneration $\cX \to \Delta$ of a smooth projective variety.\footnote{Its proof mainly consists of ``proof analysis" of the construction of the \lmhs\ in \cite{St1}, \cite{Zu}  and \cite{St2}.  Our main new point is to focus from the outset on the pair $(X,\xi)$.}  All of (i)--(iv) require knowledge of at most the \fst\ neighborhood of $X$ in $\cX$.  It is worth noting that even though the central fibre $X$ is in general not uniquely 
  definable,\footnote{Exceptions include stable curves, principally polarized abelian varieties and marked $K3$ surfaces, all of which have ``good" global moduli spaces.} the ambiguity in the \lmhs s ``washes out" in the constructions (ii), (iii), (iv).\footnote{We will not try to explain this precisely, but note that in the Clemens-Schmid exact sequence the effects of doing a modification to $X$ cancel out and leave unchanged the terms with the \lmhs s.  This phenomenon is of course familiar from Deligne's theory of \mhs s in which independence of the choice of smooth completions is established (cf.\ \cite{PS} and the references cited therein).}
  
  We   note that  the traditional approach in the study of the behavior of the \phs s in a degenerating family of smooth projective varieties is to start with a family $\cX^\ast \to\Delta^\ast$  with unipotent monodromy $T$. To this we may either   associate   a period mapping 
  \[
  \qquad \qquad\Phi:\Delta^\ast\to \Ga_T\bsl D,\qquad \qquad \Ga_T=T^\Z ,\]
  and then by \cite{Sc} to this  period mapping associate  an equivalence class of  \lmhs s. Or more algebro-geometrically we may  complete $\cX^\ast\to\Delta^\ast$ to a standard family to which by \cite{St1} we may associate the same 
  equivalence class of  \lmhs s.\footnote{More precisely, to each is associated   an  equivalence class of \lmhs s.  In \cite{St1} it is shown that the two equivalence classes of \lmhs s agree.} In this paper we are \emph{starting} with the central fibre $X$ with only the assumptions that (a) $X$ is projective and is locally a normal crossing divisor, or more generally that it is locally a product of normal crossing divisors, and (b) there exists a $\xi\in T_X\Def(X)$ that is to \fst\ smoothing and preserves the ample line bundle.  We hope that this helps to explain the title of this work.
  
    In what follows we shall use $X$ to denote both a compact analytic variety and a germ of an analytic variety; we hope the context will make clear to which we are  referring.  When $X$ is a compact analytic variety and $x\in X$ we shall denote by $X_x$ the germ  of analytic variety defined by localizing $X$ at $x$.

    The other notations we have used are either standard or will be noted where introduced.  For our variety $X$ we will have $\dim X=n$, and we shall generally consider cohomology and hypercohomology in degree $m$ (e.g.,  $H^m(X^{[k]})$'s).
    
      \section{Deformation theory}\setcounter{equation}{0}
  
  Our basic reference is \cite{Pa}, as summarized in \cite{Fr2} for the normal crossing case and whose terminology and notations we shall generally follow.\footnote{We note that the setting of log-analytic geometry is an alternate, and in many ways preferable, way to present this theory (cf.\ \cite{Abetal}, \cite{KN}, \cite{KU} and the references cited in these   works).  For example, in this context the central concept of $d$-semi-stability (cf.\ \pr{2.6} below) simply becomes the  existence of a log structure.  Moreover, a logarithmic deformation of a smoothable normal crossing variety remains smoothable; none of the ``bad" components in $\Def(X)$ can arise.  In the setting of logarithmic deformation theory, unobstructed deformations of $X$ simply means independently smoothing the connected components of $D=X_{\sing}$.  We have written this work in the traditional setting in part because this allows us more easily to connect with the other topics discussed.\label{17}}  For $X$ either a compact analytic variety, or a germ of a reduced analytic variety, we shall denote by $\Def(X)$ the  space parametrizing the corresponding family $\cX_{\Def(X)}\xri{\pi} \Def(X)$ that is versal for germs of flat families $\cX\xri{\pi} S$ with $\pi^{-1}(s_0)=X$.  The Zariski tangent space to $\Def(X)$ is 
  \[
  T_X \Def(X)= \EExt^1_{\cOX}\bp{\Om^1_X,\cOX}.\]
  As usual we think of $\xi\in T_X\Def(X)$  as giving a family
  \[
  \cxxi\to \Delta_\eps\]
  where $\Delta_\eps=\Spec \C[\eps]$, $\eps^2=0$.

  Of basic importance for us will be the  exact sequence 
  \begin{equation}\lab{2.1}
 \boxed{\begin{array}{rl} 0&\to H^1\lrp{\Ext^0_{\cOX}\bp{\Om^1_X,\cOX}}\to \EExt^1_{\cOX} \bp{\Om^1_X,\cOX}\\
 &\to H^0 \lrp{\Ext^1_\cOX\bp{\Om^1_X,\cOX}} \xri{\delta} H^2\lrp{\Ext^0_{\cOX}\bp{\Om^1_X,\cOX}}\end{array}}\end{equation}
  that results from the local to global spectral sequence for Ext.  The image of the first map will be denoted by
  \[
  T_X\Def^{\es}(X)\to T_X\Def(X);\]
  it represents the Zariski tangent space to the equisingular, or locally trivial for the germs $X_x$ in $X$, deformations.  For $x\in X$ the image of the map
  \[
  \EExt^1_\cOX\bp{\Omox}\to \Ext^1_\cOX\bp{\Omox}_x\]
  represents the $1^{\rm st}$ order deformation  of the germ $X_x$ of analytic variety induced by a global $1^{\rm st}$ order deformation of $X$.
  
  In our situation where $X$ is locally a product of normal crossing varieties given by \pr{1.2} the local deformation theory is particularly harmonious.  Taking first the case when $X$ is a germ of a normal crossing variety given locally in $\C^{n+1}$ by
  \begin{equation}\lab{2.2}
  f(x)=: x_1\cdots x_k =0\end{equation}
  with versal deformation space $\cX\subset \C^n\times \C$ given by $f(x)=t$, and 
  with the notations
  \[
  \bcs
  X_i\!=\{x_i=0\},\\
  D_i\!= X_i\cap \lrp{\bigcup_{j\ne i}X_j}&\!\hbox{defined by the ideal } \{\part_{x_i}f\!=\!x_1\cdots \wh{x_i}\cdots x_k\},\\
  D\!= \cup D_i=X_\sing\\
  X^{[\ell]}=\coprod_{|I|=\ell}X_I &\hbox{where $I=(i_1,\dots,i_\ell)$ with } 1\leqq i_1<\cdots < i_\ell \leqq n\\
  &\hbox{and }X_I=X_{i_1\cap\dots\cap X_{i_\ell}}\ecs\]
  we have as $\cOX$-modules
  \begin{equation}\lab{2.3}
  \Ext^i_{\cOX} \bp{\Omox} \cong\bcs
\ker\{T_\cX\otimes \cOX\to \cOX(X)\}& i=0,\\
  (I_X/I^2_X)^\ast\otimes \cO_D& i=1,\\
  0&i\geqq 2.\ecs\end{equation}
  This follows from the Ext-sequence arising from the exact sequence
  \begin{equation}\lab{2.4}
  0 \to
  I_X/I^2_X \to \Om^1_\cX\otimes \cOX\to \Om^1_X\to 0 
 \end{equation}
  which, setting $\vp_i= \part_{x_i}f  \, dx_i= x_1\cdots \wh{x_i}\cdots x_k\, dx_i$, gives that $\Om^1_\cX$ is  freely generated over $\cOX$, and that $\Om^1_X$ is generated  by $dx_1,\dots,dx_{n+1}$ subject to the defining relation
  \[
  df = \sum_i \vp_i =0.\]
  
  Assuming now that $X$ is a complete algebraic variety that is locally a normal crossing variety given locally by \pr{2.2}, motivated by the middle equation in \pr{2.3} and taking into account the scaling of $f$ under $f\to uf$ where $u\in\cO^\ast_X$  and following \cite{Fr2}, we may  define the \emph{infinitesimal normal bundle} by
  \begin{equation}\lab{2.5}
  \cOD (X)=\Ext^1_{\cOX}\bp{\Omox} .\end{equation}
  The point here is that, unless we are given a global embedding of $X$ as a hypersurface in a smooth variety $\cX$, we cannot define the normal bundle $\cOX(X)$, but we are able to intrinsically define what would be the restriction to $D$ of the normal bundle of $X$ in a smooth ambient space if such exists.  
  
  In more detail,  we set
  \begin{align*}
  \cO_{D_i}(-X)&= \lrp{I_{D_i}/I^2_{D_i}} \otimes_{\cO_{D_i}} \lrp{I_{X_i}/ I_{X_i} I_{D_i}}\\
  \cOD(-X)&= \lrp{I_{X_1}/I_{X_1} J_D}\otimes_{\cOD}\dots \otimes_{\cOD}\lrp{I_{X_k}/I_{X_k}J_D}\end{align*}
  where $I_{X_i}$ is the ideal sheaf of $X_i$ in $X$ and $J_{D_i}$ is the ideal sheaf of $D_i$ in $X$.  The second equation then serves to define $\cOD(X)$ in agreement with \pr{2.5} and we have \[
\Ext^1_\cOX\bp{\Omox} \cong \cO_D(-X)^\ast\otimes \cO_D\]
where $\cO_D(-X)= I_X/I^2_X$ in case we have $X\subset \cX$.

  As in \cite{Fr2}, $X$ is said to be \emph{$d$-semi-stable} if 
  \begin{equation}
  \lab{2.6}
  \cOD(X)\cong \cO_D,\end{equation}
  that is there exists a nowhere vanishing section of the line bundle $\cO_D(X)$ over $D$.\footnote{Here the point is that if we have $X\subset \cX$ with $\cX$ smooth and $\cX\xri{\pi}\Delta$ with $\pi^{-1}(0)=X$, then the conormal bundle $\cO(-X)=I_X/I^2_X$ is trivial.  Thus if we just have $X\subset \cX$ where $\cX$ is smooth, a necessary condition for there to exist $\cX\xri{\pi}\Delta$ as above is that $\cOX(X)\cong \cOX$.  The $d$-semi-stability condition \pr{1.6} is intrinsic to $X$ and does not require the existence of an $\cX$.}

  \smallbreak \centerline{\emph{We shall assume throughout that $X$ is $d$-semi-stable.}}  \smallbreak

  Returning to \pr{2.1} and using \pr{2.5} we have the map
 \vspth \[
  \EExt^1_{\cOX}\bp{\Omox}\to H^0\bp{\cOD(X)}, \vspth\]
  which we denote by $\xi\to \xi_D$ where $\xi\in\EExt^1_\cOX\bp{\Omox} = T_X\Def(X)$.  We shall say that the $1^{\rm st}$ order deformation $\xi$ of $X$ is \emph{smoothing} if $\xi_D$ is nowhere vanishing.  We shall also  generally abuse terminology by dropping the ``$1^{\rm st}$ order'', and we  shall say that $\xi_D$ is non-zero rather than nowhere vanishing.  Given such a $\xi$ we have a family $\cxxi\to \Delta_\eps$, where $\Delta_\eps=\Spec\C[\eps]$ with $\eps^2=0$, where $\cxxi$ is smooth and in which the fibre  over $0$ is $X$.  As noted earlier, it may or may not be the case that $\cxxi\to\Delta_\eps$ can be lifted to a family $\cX\to \Delta$; this issue will play no role in what follows.

  For later reference we note that the sheaf $\Om^1_X$ of K\"ahler differentials  is defined  by \pr{2.4},  where the injectivity of the first map is a property of $X$ as given by \pr{2.2}. It   is not locally free as a sheaf of $\cOX$-modules, but rather has a torsion subsheaf
  \[
  \tau^1_X\subset \Om^1_X\]
  which is locally generated by the forms $\vp_i$ define above.  Its support is $D=X_\sing$, and as noted in \cite{Fr2} since $\Om^1_X/\tau^1_X$ is locally free the above inclusion induces an isomorphism
  \[
  \Ext^1_{\cOX}\bp{\Omox}\simto \Ext^1_{\cOX}\bp{\tau^1_X,\cOX},\]
  which ``explains" the identification \pr{2.5}.
  
  Of   importance for this work will be to consider the set $A$ of connected components $D_a$, $a\in A$, of $D$.  Recalling  our blanket assumption that there exists    a $1^{\rm st}$ order smoothing deformation of $X$, we will have for each $\alpha\in A$
 \vspth  \[
  \cO_{D_\alpha}(X)\cong \cO_{D_\alpha}, \vspth\]
  where the particular isomorphism depends on the  choice of a non-zero section $\xi$ of $\cO_D(X)$.
  Thus if $\xi\in T_X \Def(X)$ with restriction  $\xi_{D_a}$ to $\cO_{D_\alpha}(X)$, we see that  \emph{along the component $D_\alpha$ of $D=X_\sing$ the deformation  of $X$ given by $\xi$ is either everywhere smoothing or equisingular.}  For each subset $B\subset A$ we set
 \vspth  \[
  T^B_X \Def(X) = \lrc{\xi\in T_X\Def(X):\xi_{D_\beta}=0\hensp{for}\beta\in B}. \vspth\]
  Then $T^B_X\Def(X)$ corresponds to the deformations that are equisingular along the $D_\beta$ for $\beta\in B$.  The extremes are
  \begin{itemize}
  \item $B=\emptyset$ corresponds to the open set $T^0_X\Def(X)$ of smoothing deformations;
  \item $B=A$ corresponds to the space $T_X\Def^{\es}(X)$ of equisingular deformations.
  \end{itemize}
  In a  way that will be explained in detail later, this gives a stratification of $T_X\Def(X)$ and leads to the definition of the cone $\sig_X$ mentioned in the introduction.
  
  \begin{exam} A simple example is when $X$ is a nodal curve.  The surjectivity of the map
  \begin{equation}
  \lab{2.7}  \EExt^1_\cOX\bp{\Omox} \to \opplus_{\alpha\in A} H^0\bp{\cO_{D_\alpha}(X)} \to 0\end{equation}
  corresponds to individually smoothing the nodes.
  \end{exam}
  
  \begin{exam}
  Suppose that $X_0$ is a singular variety with isolated singular points $p_\alpha$ given by $f_\alpha(x)=0$.  We may resolve the singularities to obtain $X$ where $D$ has connected components $D_\alpha$.  The versal deformation spaces given by  $f_\alpha(x)=t_\alpha$ for the germ of $X_\alpha$ at $p_\alpha$ and for the inverse image $X_\alpha$ of $p_\alpha$ in $X$ coincide (cf.\ \cite{Pa}).  The failure of surjectivity of the first map in 
  \[
  \EExt^1_\cOX\lrp{\Omox}\to \opplus_{\alpha\in A}H^0 \bp{\cO_{D_\alpha}(X)} \xri{\delta} H^2 \lrp{\Ext^0_\cOX\bp{\Omox}}\]
  measures the obstruction to simultaneously smoothing the $p_\alpha\in X_0$.
  \end{exam}
  
  When the $p_\alpha$ are ordinary double points the dual to the mapping $\delta$ in the sequence may be computed and leads to the conditions on the simultaneous smoothing of the nodes that may be lifted to a smoothing of $X$.\footnote{If the mapping \pr{2.7} is surjective, then
  \[
  \sig_X\otimes \R\cong \rspan_{\R>0} \{i\xi_{D_1},\dots,i\xi_{(A)}\}.\]
  The reason for the ``$i$'' in $i\xi_{D_a}$ is that if we think of $\xi_{D_a}$ as giving a tangent vector to a one parameter family then $i\xi_a$ is supposed to suggest turning around the origin --- i.e., monodromy --- in the family.}
  
  In general when $X$ may be smoothed but the connected components $D_\alpha$ may not be independently smoothed, the situation is more complicated and necessitates the blowing up of $X$.  This is the situation where in the setting of logarithmic deformation theory there are obstructions and will be discussed at another time. 
  
   The example of  K3  surfaces is discussed in \cite{Fr1}, \cite{Fr2} and \cite{KN}.  
 
  We next turn to the local case  where the germ of variety $X$ is a product
  \[
  X=X^1\times\cdots\times X^k = \prod_{\mu\in U}X^\mu\]
  of normal crossing varieties.\footnote{This case is treated in \cite{Fu1} and \cite{Fu2}.}  Letting $\pi_\mu:X\to X^\mu$ denote the projection, from the isomorphism of $\cOX$-modules
  \[
  \Om^1_X\cong\oplus \pi^\ast_\mu \Om^1_{X^\mu}\]
  leading to
  \[
  \Ext^1_\cOX \bp{\Omox}\cong \opplus_\mu \Ext^1_{\cO_{X^\mu}} \bp{\Om^1_{X^\mu}, \cO_{X^\mu}}\otimes \cOX,\]
  we may extend the local theory in the evident way.  The sequence \pr{2.4} now becomes
  \[
  \xymatrix@R=.5pt{0\ar[r] & \cO_X^\ell \ar[r] & \Om^1_\cX\otimes \cOX\ar[r] &\Om^1_X\ar[r]&0\\
  &\sidecong& \sidecong& \sidecong\\
  &\oplus  \pi^\ast_\mu \cO_{X^\mu}&\oplus \pi^\ast \Om_{\cX^\mu}\otimes \cO_{X }&\oplus  \pi^\ast_\mu \Om^1_{X_\mu}.
 }\]
  For
  \[
  \bcs
  D^\mu = X^1\times\cdots\times D^\mu\times\cdots\times X^k \\
  D= \sum D^\mu= X_\sing\ecs\]
   we have
   \[
   \Ext^1_\cOX\bp{\Omox} \cong \opplus_\mu \pi^\ast_\mu \cO_{D^\mu}.\]
   The local versal deformation space is the product of the local versal deformation spaces for the factors.
   
   As in the normal crossing case one may intrinsically define an infinitesimal normal sheaf $\cN$.  In the stratification $X_{\sing,\ell}$ of $X_\sing$ by the number of singular factors in the local product of normal crossing varieties description given by \pr{1.1}, $\cN$ is a coherent sheaf whose restriction to $X_{\sing,\ell}\bsl X_{\sing,\ell+1}$ is locally free of rank $\ell$.  The definition of $d$-semi-stability may then be extended.   This will be done in the work in progress; the practical effect of assuming $d$-semi-stability is that the  to be constructed locally defined sheaves of $\cOX$-modules $\Om^1_{\cxxi}\otimes \cOX$ and $\Om^1_{\cxxi/\Delta\eps}(\log X)\otimes \cOX$ patch together to give global sheaves over all of $X$.
   
   When we consider the global situation where $X$ is locally a product as above, we retain our standing assumption that in  the map
   \[
   \EExt^1_\cOX\bp{\Omox}\to H^0 \lrp{\Ext^1_\cOX\bp{\Omox}}\]
there is $\xi\in \EExt^1_{\cOX}\bp{\Omox}$ which   is   a smoothing deformation along each component of $X_\sing$.  This does \emph{not} mean that for each germ $X_x\subset X$ the global deformations map \emph{onto} the space of local smoothings of $X_x$.  It does mean  that there is a $\xi\in \EExt^1_\cOX\bp{\Omox}$  whose local image in each
   \[
   \Ext^1_{\cO_{X^\mu}} \lrp{\Om^1_{X^\mu},\cO_{X^\mu}}\cong \cO_{D^\mu}(X^\mu)\]
 is non-vanishing.  Then the above discussion regarding the connected components of $D$ extends and will be taken up in a future work.
   
   One significant difference   in the local situation where the number $l$ of local factors is  strictly larger than one is this:  For a 1-parameter smoothing family
   \[
   \cX_\Delta\to \Delta\]
   with tangent $\xi$, the total space $\cX_\Delta$ is \emph{singular}.  This can be seen already in the local situation
   \[
   \bcs xy=t_1\\
   uv=t_2\ecs\]
   where the disc is given by $t_1/t_2=\la\ne 0$.  Then even though the total space $\cX\to \Delta_1\times \Delta_2$ is smooth, the subvariety $\cX_\Delta\subset \cX$ is singular at the origin.

   A final comment for this section: In the study of varieties that are locally products of normal crossing varieties, the necessary multi-index notations may obscure the essential points.  Our experience has been that for normal crossings the two cases
   \[
   \bcs
   xy=0\\
   uvw=0,\ecs\]
   and for products of normal crossings the cases
   \[
   xy=0,\quad uv=0\]
   capture all the essential phenomena.  The main subtlety seems to arise when we smooth the singularity to obtain $\cX$, various exact sequences over $\cO_\cX$ fail to become exact when we restrict to $X$ by tensoring with $\cOX$ and some care must be taken in the computations to keep track of this.
   
   \section{Proofs of Theorems I and I$'$} \setcounter{equation}{0}
   This initial  discussion is mainly  local.   We begin with a germ of normal crossing variety   $X$ given by \pr{2.2}.  Given a non-zero
   \[
   \xi\in\Ext^1_\cOX\bp{\Omox}\cong \cOD(X)\]
   we denote by $\cX_\xi\xri{\pi}\Delta_\eps$ the corresponding versal family
   \[
   x_1\cdots x_k =\eps,\qquad \eps^2=0\]
   and write the   extension as
   \[
   0\to \Om^1_{\Delta_\eps}\otimes \cOX \to \Om^1_{\cX_\xi}\otimes \cOX\to\Om^1_X\to 0.\]
   Here, $\cO_{\cX_\xi}$ is locally isomorphic to $\cOX[\eps]$ and $\Om^1_{\cX_\xi}$ is the free $\cO_{X_\xi}$-module generated by $dx_1,\dots,dx_{n+1}, d\eps$ modulo the relation $d\eps = \sum^k_{i=1}\vp_i$.  Unless otherwise noted the tensor products are over $\cO_{\cxxi}$. We are    setting $\Om^1_{\Delta_\eps}=\pi^\ast \Om^1_{\Delta_\eps}$ and    are writing the sequence in this way to emphasize the scaling property with respect to $\xi$.  Note that   $\Om^1_{X_\xi}\otimes \cOX$ is the $\cOX$-module with the same set of generators and defining relation, and where in computations we set $x_1\cdots x_k=0$  but do \emph{not} set $d(x_1\cdots x_k) = \sum^k_{i=1}\vp_i=0$.
   
    We may  as usual define the free $\cO_{\cX_\xi}$-module
  $
   \Om^1_{\cX_\xi}(\log X)$ with generators $dx_1/x_1,\allowbreak\dots,dx_k/x_k, dx_{k+1},\dots,dx_{n+1},d\eps/\eps$ modulo the relation $d\eps/\eps=\sum^k_{i=1} dx_i/x_i$.  Then 
   \[ \Om^1_\cxxi (\log X) \otimes  \cOX\] is freely   generated over $\cOX$ with the same set of generators and defining relation. 
   
   We next define
   \begin{equation}\lab{3.1}
   \Om^1_{\cxxi/\Delta_\eps}(\log X)\otimes \cOX\end{equation}
   to be the $\cOX$-module with the above generators and generating relation
   \[
   \sum_i dx_i/x_i=0.\]
   We will describe this intrinsically in a moment.  Here we note the crucial point that  in the case of a global normal crossing variety $X$ fixing a nowhere zero $\xi\in \EExt^1_\cOX\bp{\Omox}$ \emph{uniquely}  locally determines  a \emph{normalized generator} $\eps-x_1\cdots x_k$ of the ideal $I_\cxxi$ of $X_\xi$: If we   have $x'_i = u_ix_i$ where $u_i\in\cO^\ast_X$, it follows from the equality 
    \vspth\[\eps-x_1\cdots x_k=\eps-x'_i \cdots x'_k
  \vspth  \]
    of normalized generators that $u=u_1\cdots u_k=1$ so that $\sum dx'_i/x'_i=\sum dx_i/x_i$.  This gives   
    \begin{quote}
     \emph{\pr{3.1} may be defined by the pair $(X,\xi)$ where $X$ is a local normal crossing variety and    $\xi\in T_X\Def(X)$   is non-vanishing.\footnote{In the setting of logarithmic geometry, to define a logarithmic structure on a normal crossing variety requires $d$-stability, and then the variety is log-smooth (\cite{St2}, \cite{KN}).}
}    \end{quote}
 
   One   small point  to notice is that the natural map
   \begin{equation}\lab{3.2}
 \pi^\ast  \Om^1_\Delta \otimes \cOX \to \pi^\ast\Om^1_\Delta(\log 0)\otimes \cOX\end{equation}
   is \emph{zero}; this is because
   \[
   d \eps\otimes 1\to \eps\lrp{\frac{d\eps}{\eps}} \otimes 1=\frac{d\eps}{\eps}\otimes \eps=0.\]
       A related point is that there is a natural map of $\cOX$-modules
   \[
   \Om^1_X\to \Om^1_{\cX_\xi/\Delta_\eps}(\log X)\otimes \cOX\]
   given on generators by $dx_i\to dx_i$, and then this map has kernel $\tau^1_X$.  Since $\Om^1_{\cxxi/\Delta_\eps}(\log X)\otimes \cOX$ is locally free we know that the above map must have kernel containing $\tau^1_X$; computation shows that equality holds.
   
   We follow the usual notations
   \begin{alignat*}{5}
   X^{[1]} &= \coprod_i X_i=\wt X &&\hbox{(normalization of $X$)}\\
   X^{[2]}&= \coprod_{i<j} X_i\cap X_j=\wt X_{\sing}&\qquad &\hbox{(normalization of $X_\sing$)}\\
   X^{[3]}&=\coprod_{i<j<k} X_i \cap X_j\cap X_k\\
   &\qquad\qquad\vdots &\end{alignat*}
   with maps
 \vspth  \[
   a_j :X^{[j]} \to X, \vspth\]
   and where we set $a_1=a:\wt X\to X$.
   \medbreak
   
   \begin{demo}\lab{3.3}
   {\bf Proposition (Basic Diagram):} \em We have
\begin{small}  \[
 \hspace*{-24pt}   \xymatrix@C=1.25pc{
   0\ar[r]&\pi^\ast\Om^1_\Delta\otimes \cOX\ar[r]&\Om^1_{\cX_\xi}\otimes \cOX\ar[r]\ar[d]& \Om^1_X\ar[d]\ar[r]&0\\
   0\ar[r]& \pi^\ast\Om^1_{\Delta_\eps}(\log 0)\otimes \cOX\ar[r] \ar[d]^{\rm Res}_\sim&\Om^1_{\cxxi} (\log X)\otimes \cOX\ar[r]\ar[d]^{\rm Res}&\ar[d]^{\rm Res}\Om^1_{\cxxi/\Delta_\eps}(\log X)\otimes \cOX\ar[r]&0\\
    0\ar[r]& \cOX\ar[r]&a_\ast \cO_{X^{[1]}}\ar[d]\ar[r]&\ker\{(a_2)_\ast \cO_{X^{[2]} }\to (a_3)_\ast \cO_{X^{[3]}}\}\ar[d]\ar[r]&0\\
    & &0&0.
      }\]\end{small}
   \end{demo}
   The right-hand map in the bottom row is the truncation of a resolution of $\cOX$ that is given in the comment immediately following the proof of this proposition. 
   
   \begin{proof}
   We begin with the standard diagram of $\cO_{\cxxi}$-modules
   \begin{small} \[
   \xymatrix@C=1.25pc{
&0\ar[d]&0\ar[d]&0\ar[d]&\\
0\ar[r]&\pi^\ast\Om^1_\Delta\ar[d]\ar[r]& \Om^1_{\cxxi}\ar[d]\ar[r]& \Om^1_{\cxxi/\Delta}\ar[r]\ar[d]&0\\
0\ar[r]&\pi^\ast \Om^1(\log 0)\ar[d]^{\Res}\ar[r]& \Om^1_{\cxxi}(\log X)\ar[r]\ar[d]^{\Res}& \Om^1_{\cxxi/\Delta}(\log X)\ar[d]^{\Res}\ar[r]&0\\
0\ar[r]&\cOX\ar[d]\ar[r]& a_\ast\cO_{X^{[1]}}\ar[d]\ar[r]&\ker \lrc{ (a_2)_\ast \cO_{X^{[2]}}\to (a_3)_\ast \cO_{X^{[3]}}}\ar[d]\ar[r]&0\\
&0&0&0&}\]
   \end{small}%
   
   \noindent where the $\cOX$-modules on the bottom row are considered as $\cO_{\cxxi}$-modules by the restriction map $\cO_{X_\xi}\to \cOX$. When we restrict to $X$ by    \[
   \cF\to \cF\mid_X=\cF\otimes_{\cO_\cX}\cOX \]
   for an $\cO_{\cX}$-module $\cF$, 
then as noted above we may lose exactness in certain places.  Calculations in local coordinates gives the exactness in the basic diagram, where we note that
   \[
   \Om^1_{\cxxi/\Delta_\eps}\otimes \cOX=\Om^1_X\]
   as the top row reduces to the defining relation  $0\to I_X/I^2_X \xri{d} \Om^1_\cxxi\otimes \cOX\to \Om^1_X\to 0$ for   K\"ahler differentials. 
      
   To give the flavor of the calculations we consider the simplest non-trivial case of $xy=\eps$.  Then
   \begin{itemize}
   \item $\Om^1_\cxxi\otimes \cOX$ is generated as an $\cOX$-module by $dx,dy,d\eps$ with the defining relation $xdy+ydx=d\eps$;
   \item $\Om^1_\cxxi(\log X)\otimes \cOX$ is generated by $dx/x$, $dy/y$, $d\eps/\eps$ with the defining relation $dx/x+dy/y=d\eps/\eps$.
   \end{itemize}
   Any $\om\in \Om^1_\cxxi \otimes \cOX$ is of the form   $f(x,y)dx+g(x,y) dy$, and using that $\otimes\cOX$ means setting ``$xy=0$''  we see that $\om$ may be normalized to be
   \[
   \om = \bp{f_1(x)+f_2(y)} dx+\bp{g_1(y)+g_2(y)}dy,\quad f_1(0)=g_2(0)=0.\]
   Similarly, $\vp\in \Om^1_\cxxi(\log X)\otimes \cOX$ may be normalized to be
   \[
   \vp=\bp{a_1(x)+a_2(y)} \frac{dx}{x} +\bp{b_1(x)+b_2(y)}\frac{dy}{y},\qquad a_1(0)=b_2(0)=0.\]
   Then
   \[
   \Res \vp = a_2(y)\oplus b_1(x)\in (a_\ast)\cO_{X^{[1]}}.\]
   If $\Res\vp=0$, then writing $a_1(x)=x\wt{a_1}(x)$ and $b_2(y)=y\wt{b_2}(y)$ we have
   \[
   \vp=\wt{a_1}(x)dx + \wt{b_2}(y)dy\in\Om^1_\cxxi\otimes \cOX.\]
   
   A similar calculation gives the exactness of the right-hand column.  For the case of a triple point $xyz=\eps$ the residue calculation is more complicated and is similar to (2.10) in \cite{Fr2}.\end{proof}

  We want to make two  comments on the basic diagram.  The first is 
{\setbox0\hbox{(1)}\leftmargini=\wd0 \advance\leftmargini\labelsep
  \begin{quote}
  \em The bottom row in the basic diagram is the truncation of the resolution
  \vspth\[
  0\to\cOX\to (a_1)_\ast \cO_{X^{[1]}}\to (a_2)_\ast \cO_{X^{[2]}} \to (a_3)_\ast \cO_{X^{[3]}}\to\cdots  \vspth\]
  of $\cOX$.\end{quote} }

  \noindent 
  The map $(a_k)_\ast \cO_{X^{[k]}} \to (a_{k+1})_\ast \cO_{X^{[k+1]}}$ is given by
  \[
 u_{i_1\cdots i_k}  \hensp{on}X_{i_1}\cap\dots \cap  X_{i_k}  \]
  maps to
  \[
  \sum_j   (-1)^{j+1} u_{i_1\cdots \hat i_j\cdots i_{k+1}}\big|_{X_{i_1}\cap \cdots \cap X_{i_{k+1}}}.\]  
    This is standard (cf.\ \cite{Fr2}, \cite{St1} and \cite{Zu}). We note also  the resolution
  \vspth  \[
  0\to \C_X\to (a_1)_\ast \C_{X^{[1]}}\to (a_2)_\ast \C_{X^{[2]}}\to (a_3)_\ast \C_{X^{[3]}}\to\cdots  \vspth\]
  of the constant sheaf on $X$.
  
  The second is that  we list   the main take-aways from the basic diagram:
  \begin{demo}\lab{new3.4}
{\setbox0\hbox{(1)}\leftmargini=\wd0 \advance\leftmargini\labelsep
   \begin{enumerate}
  \item Given $X$ and $\xi\in \EExt^1_\cOX\bp{\Omox}$ with the property that $\xi$ is non-zero along $D$, we may by definition construct an extension of $\cOX$-modules
  \[
  0\to\cOX\to \cF_\xi\to \Om^1_X\to 0.\]
  \item From \cite{Pa}, we may actually construct a space $\cxxi$ with structure sheaf $\cO_\cxxi$ locally isomorphic to $\cOX[\eps]$ giving a mapping $\cxxi\to \Delta_\eps=\Spec \C[\eps]$ where
  \begin{itemize}
  \item $\cxxi$ is smooth (this is the assumption that $\xi_D\ne 0$);
  \item $\cF_\xi\cong \Om^1_\cxxi\otimes \cOX$ as $\cOX$-modules; this is the top row in the basic diagram.\end{itemize}
  \item We may then proceed, using $\Om^1_\cxxi$ as an $\cO_{\cxxi}$-module, to construct the remainder of the basic diagram; the inclusion map  $\cOX\to \Om^1_{\cX_\xi} (\log X)\otimes \cOX$ is given by $1\to d\eps/\eps$, and then the quotient defines  the $\cOX$-module $\Om^1_{\cxxi/\Delta_\eps}(\log X)\otimes \cOX$.
  \item From this we may, in the standard way,    proceed to construct the complex $\bp{\Om^\bullet_{\cxxi/\Delta_\eps}(\log X)\otimes \cOX,d}$; as will be noted below, the hypercohomology of this complex will give the complex vector space  $V_\xi$ and Hodge filtration $\fbul_\xi$  for the \lmhs.
  \item The previous steps are  either explicit or implicit in \cite{Fr2}; the final  steps to define the weight filtration and $\Q$-structure may then be carried out by  the methods in \cite{St2}.\end{enumerate}}\end{demo}
  \smallbreak
  \noindent We will elaborate more on this at the end of the section.
  
   Turning now to the case where $X$ is locally a product of normal crossing varieties as given by \pr{new1.2}, we may extend the discussion above with one significant change.  Namely, in the local situation instead of a single smoothing deformation $\xi\in T^0_X\Def(X)$ we now need to be given a  $k$-tuple
  \[
\bxi= (\xi_1,\dots,\xi_k)\]
  where $\xi_i$ smooths the factor $X_i$ in $X$.  Then for $\la=(\la_1,\dots,\la_k)$ with all $\la_i\ne 0$ 
  \[
  \bxi_\la =:\sum \la_i \xi_i \in T_X\Def(X)\]
  is a smoothing deformation of $X$.  This is all local. 
  
  Globally we need to be given an $\ell$-tuple $\bxi \in\opplus^\ell T_X\Def(X)$ such that locally around each $x\in X$ there is a  $k$ sub-tuple of $\bxi  $ that satisfies the above condition.
  
    We note again the difference when the number of local factors $k\geqq 2$; if $\bxi_\la$ is tangent to a family $\cX_{\Delta_\la}\to \Delta_\la$, then
  \begin{equation}\lab{3.4}
  \hbox{\em the total space $\cX_{\Delta_\la}$ is singular.}\end{equation}
  These singularities are of a standard form and may be resolved to give a standard family 
  \[
  \wt{\cX}_{\Delta_\la}\to \wt \Delta_\la \] 
  where $\wt \cX_{\Delta_\la}$ is smooth.   
    
  Setting $\Delta_{\eps_j }= \Spec \C[\eps_j]$ and $\Delta_\beps=\prod_j \Delta_{\eps_j }$, using the projection $X_1\times \cdots \times X_k\xri{\pi_j} X_j$, we define
  \begin{equation}\lab{3.5}
  \Om^1_{\cxxib/\Delta_\beps} (\log X)\otimes \cOX =\oplus \pi^\ast_j\Om^1_{\cX_j/\Delta_{\eps_j}} (\log X_j) \otimes \cO_{X }.\end{equation}
  In coordinates, for the case where $X$ is given by
  \[
  \bcs xy=t_1\\
  uv=t_2\ecs\]
  so that $I_X$ is generated by $xy$ and $uv$,
   $\Om^1_{\cxxib/\Delta_\beps}(\log X)\otimes \cOX$ is generated as an $\cOX$-module by $dx/x$, $dy/y$, $du/u$, $dv/v$ with the relations $dx/x+dy/y=0$, $du/u+dv/v=0$.  This coordinate description  extends in the evident way when $X$ is given by \pr{new1.2}.
   
   Finally, we will relate this construction to that given in \cite{Fu1}, \cite{Fu2}.  We have 
   \begin{equation}\lab{3.6}
   \cX\xri{\pi} S\end{equation}
   where   locally in $\C^{n+k}$ with coordinates $(x_1,\dots x_n, t_1,\cdots , t_k)$ and using the notation \pr{new1.2}, $\cX$ is given by
   \begin{equation}
   \lab{3.7}
   \bcs x_{I_1}&\chs=t_1
   \\
& \chs\vdots\\
   x_{I_k}&\chs = t_k\ecs\end{equation}
   and $\pi$ is the projection $(\bx,\bt)\to \bt$.  There are then normal crossing divisors $\cY\subset \cX$ and $T\subset S$ such that \pr{3.6} is a map
   \[
   (\cX,\cY)\to (S,T)\]
   as defined in \cite{Fu1}, \cite{Fu2}.  If $\dim X=n$ and $\dim S=\ell$, then $\dim \cX=n+\ell$.  Locally $S$ is embedded in $\C^k\times \C^{\ell-k}$ where the first $k$ coordinates are the $t_i$ above  and the remaining $\ell-k$ coordinates are parameters.   We note that $\pi^{-1}$ (set of coordinate hyperplanes in $\C^k$) is a singular subvariety of $\cX$.  Globally, we will have divisors $D_1,\dots,D_\ell$ on $\cX$ such that locally $D_1,\dots D_k$ are the inverse image under $\pi$ of the coordinate hyperplanes  $t_i=0$ and $\cY=D_1+\dots +D_\ell$ is a reduced normal crossing divisor in $\cX$ with $X=D_1\cap\cdots \cap D_\ell$.
 
   \subsubsection*{Discussion of the proofs of Theorems {\rm \ref{thm1}} and {\rm \ref{thm1p}}}\hspace*{-6pt}\footnote{In a work in progress we intend to provide details for this argument with emphasis on the local structure and how this relates to the results in \cite{CK} and \cite{CKS1}.}
   For the case when $X$ is a local normal crossing variety, using \pr{new3.4}   the essentials of the proof are in \cite{Fr2} and  \cite{St2}.  The sheaves
   \[
   \wedge^\bullet_\xi = \wedge^\bullet \Om^1_{\cxxi/\Delta_\eps}(\log X)\otimes \cOX\]
   form a filtered complex in the evident way, and
 \[
   \bcs
   V_\xi=  \bH^m(\wedge^\bullet_\xi)\\
   \fbul V_\xi= \fbul \bH^m(\wedge^\bullet_\xi)\ecs\]
    defines the vector space and Hodge filtration for the \lmhs.  As usual, dating to \cite{St1} (cf.\ also \cite{Zu}), the construction of the monodromy weight filtration and $\Q$-structure are more subtle.  These may be  carried out by an adaptation of the methods  in \S 5 in \cite{St2}.
    
    More specifically, in \cite{St1} and \cite{Zu} associated to a standard family $\cX\to \Delta$ several cohomological mixed Hodge complexes are constructed.  One of these, denoted there by $\abul$ (recalled in the proof of Theorem \pr{thm7} in Section IV below) leads to the limit \mhs. Another of these, denoted by $L^\bullet$ in loc.\ cit., leads to the \mhs s on $H^\ast(X\bsl\cX)$ and on $H^\ast(\cX,\cX\bsl X)$.  In \cite{St2}, in the setting of log geometry which in his Section 5  corresponds to our $(X,\xi)$, the analogue of $L^\bullet$, denoted there by $K^\bullet$, is constructed.  Analysis of the construction leads to a cohomological mixed complex in our $(X,\xi)$ setting that gives a \lmhs\ on $V_\xi$.
   
   For the general case  where  the central fibre  $X$ in a global map \pr{3.6} is  given locally by \pr{3.7}, in \cite{Fu1}, \cite{Fu2} the methods of \cite{St1} are extended to show that for $s\in S$ and $X_s=\pi^{-1}(s)$ the hypercohomology of the complexes $\Om^\bullet_{\cX/S}(\log \cY)\otimes \cO_{X_s}$ give \mhs s.  The adaptation of the calculations there extending the methods in \cite{St2}   to the several variable log-geometry setting that corresponds to our situation will then give  the result.  As we have no substantive content to add to what is implicit  in \cite{Fu1}, \cite{Fu2} and \cite{St2} here we will not  write out the details, but rather defer them to a later work.
   
   An outstanding issue, as noted in \cite{Fu1}, \cite{Fu2}, for a family \pr{3.6}, is  in what way the \mhs\ constructed in \cite{Fu1}, \cite{Fu2} using  $\Om^\bullet_{\cX/S}(\log X)\otimes \cOX$ relates to the \lmhs\ given along a disc $\Delta_\la$ in $S=\Delta_1\times\dots \times \Delta_\ell$ in \cite{CKS1}.   One main point may be \pr{3.4}. In case $X$ is a local normal crossing variety, we have noted that $\cX_{\Delta_\bla}$ is smooth and $\cX_{\Delta_\bla}\to \Delta_\bla$ is a standard family, so the result that the \lmhs s are the same is true in this case.  Another outstanding matter is the construction of the monodromy logarithms $N_i$ from the $dt_i/t_i$ in the complexes constructed in \cite{Fu1}, \cite{Fu2}, and then to show that these give the structure    as in \cite{CKS1}.  This also will be taken up in a later work.

   Another issue, one that arises already when $X$ is a normal crossing variety whose singular locus $D=\cup D_\alpha$ has connected components $D_\alpha$, is this:  In the exact sequence \pr{2.1} when the mapping
   \[
   \oplus H^0(\cO_{D_\alpha}) \xri{\delta} H^2\lrp{\Ext^0_{\cOX}\bp{\Om^1_X,\cOX}}\]
   is non-zero, the Kuranishi space may be unobstructed but the $D_\alpha$ cannot be individually smoothed.\footnote{In the logarithmic deformation theoretic context, there are non-zero obstructions in the logarithmic analogue of $T_X\Def(X)$.}  Suppose for example that there are three components so that projectively
   \[
   \P \lrp{\oplus H^0 \bp{\cO_{D_\alpha}}} =\P^2\]
   pictured as
   \[
   \begin{picture}(60,60)
   \put(0,10){\line(1,0){60}}
   \put(10,0){\line(0,1){54}}
   \put(0,50){\line(4,-3){64}}
   \put(10,10){\circle*{3}}
   \put(54,10){\circle*{3}}
   \put(-3,35){$\scriptstyle  D_3$}
   \put(11,2){$\scriptstyle D_1$}
   \put(45,2){$\scriptstyle D_2$}
   \put(10,42){\circle*{3}}\end{picture}\]
   where the vertices correspond to the $H^0\bp{\cO_{D_\alpha}}$.  If $\dim (\ker \delta)=2$, there are the following possibilities for the dotted line $L=\P(\ker\delta)$:
 \[\begin{picture}(355,250)
 \put(0,220){(i)}
 \put(0,120){(ii)}
 \put(0,20){(iii)}
 \put(50,0){$\includegraphics{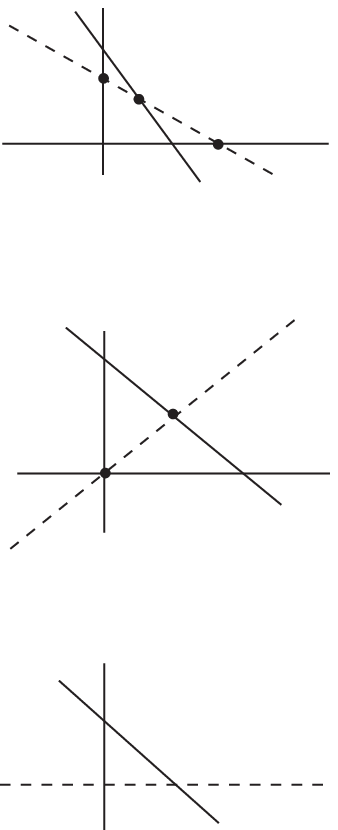}$}
 \end{picture}
   \]
   Here, (iii) does not occur since we are assuming that $X$ may be smoothed.  For (ii), assuming that $T_X \Def(X)$ is unobstructed we have a 2-parameter family $\cX\to \Delta\times \Delta$ where along one axis $D_1$ is smoothed while $D_2$ and $D_3$ deform equisingularly.  Along the other axis a similar thing happens with the roles of $D_1$ and $D_2,D_3$ interchanged.
   
   In case (i) we have a 2-parameter family with three axes along each of which one pair from $D_1,D_2,D_3$ deforms equisingularly while the remaining component of  $D$ is smoothed.  Thus the picture of the tangent space to the 2-parameter family  $\cX\to S$ is
   \[
   \begin{picture}(100,75)
   \put(0,0){\includegraphics{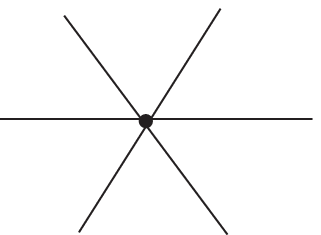}}
   \put(100,35){,}\end{picture}\]
   and for  the family $\cX^\ast \to S^\ast \cong \Delta\times \Delta\bsl \{\hbox{3 lines}\}$ where the fibres are smooth we have
   \[
   \pi_1(S^\ast)\cong\pi_1(\P^1\bsl \{0,1,\infty\}).\]
   Thus to arrive at a CKS situation we have to blow up $S$ at the origin and resolve singularities to arrive at a standard situation $\wt X\to \wt S$ where locally $\wt S$ is a $\Delta\times \Delta$ with singular fibres over the two axes.  The fibre $\wt X$ over the origin is related to $X$ in a standard way; it will have as one component a branched covering of the desingularization of $X$ along two of the $D_\alpha$, and the other components are easily described (if $X_3=\emptyset$ they are the projectivized normal bundles of the $D_\alpha$ in $\cX$).  We are now back in the situation of CKS but with a different $X$.  As in \cite{KN}  for interesting examples this complexity does not arise, and for theoretical purposes we can at least begin by assuming that $\delta=0$ in \pr{2.1}. 
          \section{Proof of Theorem II}   \setcounter{equation}{0}
      
      We first will consider the question,  informally stated as
\[
\hbox{\em What do we mean by $T$(LMHS)?}\]
Here, ``LMHS'' is the set of \lmhs s with monodromy $N$.

Setting 
\[
D_N= D\cup B(N),\]
in \cite{KU} there is  defined on $D_N$ the structure of a ``log-analytic varity with slits."  In particular, the tangent space $T_{[\fbul]} D_N$ to the underlying analytic variety at a point $[V,\wbul(N),\fbul]$ is defined, where the brackets denote the equivalence class of nilpotent orbits under the equivalence  relation $\fbul\sim \exp(zN)\cdot \fbul$.  This information may be refined if we do not pass to equivalence classes.  Thus we define
\[
\wt D_N=D\cup \wt B(N)\]
and seek to define $T_{\fbul}\wt D_N$.  For the subspace $T_\fbul \wt B(N)$ of $T_\fbul\wt{D}_N$, recalling  that $\wt D_N\subset \CD$ we set 
\begin{equation}\lab{4.1}
T_{\fbul}\wt D_N=\{ \tau\in T_{\fbul} \CD:\tau(F^p)\subset F^{p-1}\} = T_{\fbul} \CD.\footnotemark\end{equation}
\footnotetext{Essentially we are interpreting the additional infinitesimal information that is present if we consider $D_N$ as a log-analytic variety.}

We next consider the question
\[
\hbox{\em What is the algebro-geometric analogue of \pr{4.1}?}\]
This means: What algebro-geometric object maps to \pr{4.1}, extending what  is given for a smooth $X$ by \pr{1.1}?

For this we recall that associated to a pair $(X,\xi)$, where $X$ is locally a product of normal crossing varieties and $\xi\in T^0_X \Def(X)$, is a \lmhs\ whose underlying vector space is $\bH^m(\Obul_{\cxxi/\Delta_\eps} (\log X)\otimes \cOX)$.  With the identification $T_X \Def(X) = \EExt^1_\cOX(\Omox)$ in mind, we \emph{define} 
\begin{equation}\lab{4.2}
T_{(X,\xi)}\Def(X,\xi) = \EExt^1_\cOX\bp{\Om^1_{\cX_\xi/\Delta_\eps}(\log X)\otimes \cOX,\cOX}\end{equation}
where $\Def(X,\xi)$ is the set of deformations of the pair $(X,\xi)$.\footnote{As was noted in the introduction, the right-hand side of \pr{4.2} appears naturally in logarithmic deformation theory.}  We shall not attempt here to give a precise definition of $\Def(X,\xi)$, but rather shall simply take \pr{4.2} as the definition of its  tangent space.  As partial justification, we observe that with this definition the obvious map
\[
\Def(X,\xi)\to \Def(X),\]
 together with the map $\Om^1_X\to \Om^1_{\cX_\xi/\Delta_\eps}(\log X)$, give
 \begin{equation} \lab{4.3}
 \begin{array}{ccc}
 T_{(X,\xi)}\Def(X,\xi)&\longrightarrow & T_X \Def(X)\\
 \sideeq&&\sideeq\\
 \EExt^1_{\cOX}\lrp{\Om^1_{\cxxi/\Delta_\eps}(\log X)\otimes \cOX,\cOX}&\longrightarrow&\EExt^1_{\cOX}\bp{\Omox}.\end{array}\end{equation}

Turning to the definition of the maps in Theorem \ref{thm2}, from the middle row in the basic diagram \pr{3.3} one may in the usual way infer the exact sequence of complexes
\[
0\to \Om^{\bullet-1}_{\cxxi/\Delta_\eps}(\log X)\otimes \cOX \to \Obul_{\cxxi} (\log X)\otimes \cOX\to \Obul_{\cxxi/\Delta_\eps}(\log X)\otimes \cOX\to 0.\]
The connecting homomorphism in the long exact hypercohomology sequence induces 
\[
\bH^m \lrp{\Obul_{\cxxi/\Delta_\eps}(\log X)\otimes \cOX}\xri{\nabla_\xi} \bH^{m+1}
\lrp{\Om^{\bullet-1}_{\cxxi/\Delta_\eps} (\log X)\otimes \cOX}\]
which satisfies
\[
\nabla_\xi  F^p_\xi \subset F^{p-1}_\xi.\]
Then from \pr{1.5} the element $\xi^{(1)} \in T_{(X,\xi)}\Def(X,\xi)$ gives an extension class in the above exact sequence of complexes, and using the identification
\[
\bH^{m+1}\lrp{\Om^{\bullet-1}_{\cxxi/\Delta_\eps}(\log X)\otimes \cOX} = \bH^m\lrp{ \Obul_{\cxxi/\Delta_\eps} (\log X)\otimes \cOX}\]
we obtain the map in the statement of Theorem \ref{thm2}.  The fact that we map to $\End_{\rm LMHS}$ is a consequence  of the naturality of the construction of the \lmhs .

The geometric picture to keep in mind is   this: The kernel of the map \pr{2.1}
\[
\EExt^1_\cOX \bp{\Omox}\to H^0 \Bp{\Ext^1_\cOX\bp{\Omox}}\]
represents the tangents to the equisingular deformations of $X$.\footnote{Recall that we are assuming that under any non-smoothing deformation $X'$ of $X$, including an equisingular one, the deformed $X'$ remains smoothable; the condition for this is in \cite{Fr2}.}  Modulo this kernel, the image of the above map reflects how the singularities are deforming.  Given a smoothing deformation $\xi$, we may think of $\xi^{(1)}$ as  giving us the infinitesimal change in this picture.\footnote{Of course, there is more information than this in $\xi$.}  We will now illustrate this by example where it will be quite clear how  the map in the statement of Theorem \ref{thm2} gives information beyond that in the   differential 
\[
T_{s_0}S\to T_{[\fbul]}   D_N.\]
  
   Here we are imagining a family $\cX\to S$ where $S=\Delta^\ell$ and where the fibres are smooth over $S^\ast = \Delta^{\ast \ell}$ with commuting monodromy logarithm transformations $N_1,\dots ,N_\ell$ around the axes.  The corresponding nilpotent orbit  is $\exp(z_1N_1+\cdots z_\ell N_\ell)\cdot  \fbul$.
 
\begin{exam}
This will be a simpler version of the  example from the beginning of Section VI, and we  will use the notations from there.  Then $\fbul$ is a single $F$ given by the span  of the columns in the matrices below:
{\allowdisplaybreaks\begin{align*}
F&\longleftrightarrow \bpm 1\\
&1\\
&&1\\
a_{11}&a_{12}&b_1\\
a_{21}&a_{22}&b_2\\
b_1&b_2&c\epm ,\qquad a_{12}=a_{21}\\
[F]&\longleftrightarrow \bpm 1\\
&1\\
&&1\\
0&a_{12}&b_1\\
a_{21}&0&b_2\\
b_1&b_2&c\epm\end{align*}
where the notation $\longleftrightarrow$ means ``corresponds to.''  For $[F]$ we have normalized the point on the several variable nilpotent orbit by $a_{11}=a_{22}=0$:
\begin{align} \lab{new4.4}
T_{F}\CD&\longleftrightarrow
\bpm da_{11}&da_{12}&db_1\\
da_{21}&da_{22}&db_2\\
db_1&db_2&dc\epm,\\ \lab{new4.5}
T_{[F]} D_N&\longleftrightarrow \bpm 0&da_{12}&db_1\\
da_{21}&0&db_2\\
db_1&db_2&dc\epm.\end{align}}%
Geometrically, the second contains the information in $T_X \cC\subset T_X\part \cM_3$, while the first contains this information plus the information in the normal  space to $\cC$ in $\ol{\cM}_3$; i.e., the refined direction  of approach to $X$ in the boundary of $\cM_3$.  Here the term ``refined direction of approach" means the following: The crude normal direction of approach to $X$ is given by $N_1,N_2$, which may be thought of as the normal direction of approach to the image of $X$ in $D_N$.  The refined direction of approach picks out more subtle information beyond that given simply by the logarithmic terms in the period matrix.

Still referring  to the next section for the notations, the $N$-strings associated to the \lmhs\ in this example may be written as
\[
\begin{array}{c}
H^0(-1)\to H^0\\[4pt]
H^1.\end{array}\]
\end{exam}
\noindent The extension data in $\Ext^1_{\rm MHS}(\Gr_1,G_0)=\Ext^1_{\rm MHS}\bp{H^1(\wt X),H^0(D)}$ corresponds to $b_1,b_2$, while that for $\Ext^1_{\rm MHS}\bp{H^0(D)(-1),H^0(D)}$ corresponds to the $2\times 2$ symmetric matrix $(a_{ij})$.  In this case only the off-diagonal terms are invariant under $\fbul\to \exp(z_1 N_1+z_2N_2)\cdot\fbul$, while the diagonal terms require the choice of $\xi$.  

 From a cohomological perspective, the $F$ in the \lmhs\ is
 \[
 H^0 \lrp{\Om^1_{\wt X}\bp{\log (p+q)}}\]
 and the matrix $dF$ in \pr{new4.4} is in
 \[
 \Hom_s (F,\C^4/F) \cong \Hom_s\lrp{H^0\bp{\Om^1_{\wt X}(\log (p+q)), H^1(\cOX)}}\footnotemark\]
 \footnotetext{Here we are identifying $H^1(\cOX)$ with $V_{\xi/F}$ where $(V_\xi,\wbul(N),F)$ is the \lmhs\ with $N=N_1+N_2$.}
 \noindent 
 where $\Hom_s$ are the symmetric maps.
 Under the inclusion
 \[
 H^0(\Om^1_{\wt X}) \hookrightarrow H^0 \lrp{\Om^1_X(\log (p+q))}\]
 the matrix \pr{new4.5} contains a part in
 $\Hom \lrp{H^0(\Om^1_{\wt X}),H^1 (\cOX)}$.  The term $dc$ is in 
 \[\Hom\lrp{H^0\lrp{\Om^1_{\wt X}}, H^1(\cOXt)}\] and $(db_1,db_2)$ belongs to 
 \[\Hom \lrp{H^0\lrp{\Om^1_{\wt X}}, H^1( \cOX)/H^1\lrp{\cO_{\wt X}}} .\]
 
 An extreme example of the extra information is given by the genus 2 curve degenerations
 \[
 \begin{picture}(330,240)
 \put(20,0){\includegraphics{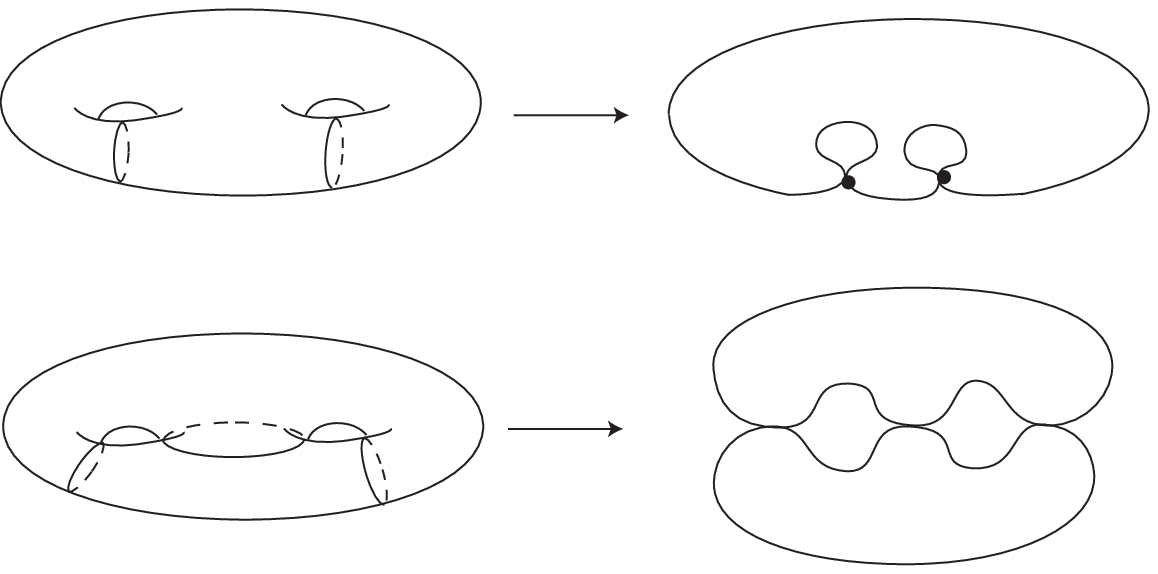}}
 \put(-10,183){(i)}
 \put(48,145){$\delta_1$}
 \put(115,145){$\delta_2$}
 \put(40,60){$\delta_1$}
 \put(79,72){$\delta_3$}
 \put(130,60){$\delta_2$}
 \put(-10,93){(ii)}
 \end{picture}\vspace*{-34pt}\]
 In each case the polarized \lmhs\ is
\[
\opplus^2 \Q(-1)\to \opplus^2 \Q\]
and the $\Ext^1_{\rm PLMS}(\bbul)$ is given by a $2\times 2$ symmetric case.  For $T_{\fbul}\wt D_N$ and $T_{[\fbul]}D_N$ they are
\[
\begin{matrix}
T_{\fbul}\CD_N\\[24pt]
T_{[\fbul]}D_N\\[30pt]\end{matrix}\qquad\qquad
\begin{matrix}
\bpm a&b\\ b&a\epm\\ \\[-4pt]
\bpm 0&b\\
b&0\epm\\ \\[-4pt]
{\rm (i)}
\end{matrix}
\qquad \qquad
\begin{matrix}
\bpm a&b\\
b&a\epm\\ \\[-4pt]
\bpm 0&0\\
0&0\epm\\ \\[-4pt]
{\rm (ii)}\end{matrix}\]
This illustrates the additional information contained in considering the map $T_{(X,\xi)} \Def(X,\xi)\to T_{\fbul}\wt D_N$.
 
 The above examples are  of course   special.  However, the regularization of logarithmic integrals  phenomenon they illustrate are fairly general.  For instance, in \cite{GGR} the generic degenerations  of Hodge structures of odd weight $n=2m+1$ are given, for $1\leqq k\leqq n$, by a specialization
 \[
 X_t \to X_0\]
 where locally in $\C^{2m+2}$   $X_0$ is  given  by
 \[
 x_1 x_2+\cdots + x_{2k-1} x_{2k}=0 \]
 and thus has a double locus of codimension $k$. 
 For $k=1$ we have that $X_0 = X$ has a codimension one  double $2m$-fold $X_\sing = X_1$. For  $k=m$, $X_0$ has an ordinary isolated quadratic singularity.  For $k\geqq 2$ we have to blow up $X_0$ to achieve a standard family.  The \lmhs s are 
 \begin{align*}
 k=m\qquad &\left\{ \begin{matrix}
 H^0(-m-1)&\longrightarrow & H^0(-m)&\qquad \dim H^0=1\\[6pt]
& \vdots  \\
 & {H^n_m}&& \\[6pt]
 \end{matrix}\right.\\
 k=1\qquad&\left\{ \begin{matrix}
 \enspace H^{2m-1}(-1)&\longrightarrow & H^{2m}\\[6pt]
 &\vdots
 \\
 & {H^n_1}\end{matrix}\right.\end{align*}
 For $ {H^n_m}$ the Hodge numbers  $h^{p,q}_k$ are the same as the $h^{p,q}$ for the original \phs\ on $ {H^n}(X_t)$, except that
 \[
 h^{m+k,m-k+1}_k = h^{m+k,m-k+1}-1.\]
 Geometrically we have a class $\om_t$ in $H^0(\Om^n_{X_t})$ that acquires a pole of  order $k$ along $X_\sing$, and by a residue-type construction we end up with a class in $H^{n-m}\lrp{\Om^{m-k+1}_{X^{[m]}}}$. The above analysis of regularizing an integral
 \[
 \lim_{t\to 0} \int_{\ga_t}\om_t\]
 then will carry over.  The details of this   will be carried out in a future work.
 
Finally we would like to give a general  cohomological description of the extra information in the map in Theorem \ref{thm2}.  We will do this in case $X$ is a nodal curve; this description will extend to the general case when $X$ has only an ordinary double locus $D=X_\sing$.  With this assumption  the right-hand column in the basic diagram \pr{3.3} gives a map
\[
\Om^1_{\cxxi/\Delta_\eps}(\log X)\otimes \cOX\to \cOD,\]
which induces a map
\begin{equation}
\lab{4.6}
\EExt^1_\cOX\bp{\cOD,\cOX}\to \EExt^1_\cOX\lrp{\Om^1_{\cxxi/\Delta_\eps} (\log X)\otimes \cOX,\cOX}.\end{equation}
The   image of this map represents the ``extra information" contribution to the map in Theorem \ref{thm2}.  The local result we need to describe this for $X$ a germ given by $xy=0$ is 
\begin{gather}\lab{4.7}
\Ext^1_{\cOX}\lrp{\cOD,\cOX}\cong \cOD\hbox{ \em and is generated by the extension}\\
\notag
0\to \cOX\to (a_1)_\ast \cO_{X^{[1]}} \to \cOD\to 0.\end{gather}
In the case under consideration,
\[
\EExt^1_\cOX\bp{\cOD,\cOX}\cong H^0\lrp{\Ext^1_\cOX\bp{\cOD,\cOX}}\]
 and there are  the number   of connected components of $D$ additional parameters picked up in the additional information.
 
 A cohomological formulation that identifies the $N$ in a \lmhs\ is this.  Recall the bottom two rows in the basic diagram \pr{3.3}, where if we use the notation $\cS = \ker\{(a_2)_\ast \cO_{X^{[2]}}\to (a_3)_\ast \cO_{X^{[3]}}\}$ and identify $\pi^\ast \Om^1_{\Delta_\eps}(\log 0)\otimes \cOX$ with $\cOX$ we have
 \begin{equation}\begin{array}{l}\lab{4.8}
 \xymatrix{0\ar[r]&\cOX\ar[r]\ar@{=}[d]& \Om^1_{\cOX}(\log X)\otimes \cOX\ar[r] \ar[d]^{\Res}&\Om^1_{\cX/\Delta_\eps}(\log X) \otimes \cOX\ar[r] \ar[d]^{\Res}&0\\
 0\ar[r]& \cOX\ar[r]&(a)_\ast \cO_{X^{[1]}}\ar[r]& \cS\ar[r]&0.}
 \end{array}\end{equation}
 This gives
 \[
 \EExt^1_{\cOX} (\cS,\cOX)\xri{\Res^\ast} \EExt^1_{\cOX} \lrp{\Om^1_{ \cX/\Delta_\eps}(\log X)\otimes \cOX,\cOX}
  \to F^{-1}\End_{\rm LMHS}.\]
  Then by interpreting the construction in \cite{St1} we find that
  \begin{quote}
  \em the image of the extension class in the bottom row of \pr{4.8} is the monodromy logarithm $N$.\end{quote}
  
Finally, we would like to point out the paper \cite{Ca-Fe} in which the notion of an infinitesimal variation of Hodge structure at infinity is defined.  Their definition pertains to equivalence classes of  \lmhs s for several variable nilpotent orbits as in \cite{CKS1}.  In the above example the  definition in \cite{Ca-Fe} would record the data
\[
\{ N_1,N_2:d a_{12},\, db_1,db_2,dc\}.\]
Roughly speaking, this data corresponds to $TB(N)$ and to the normal space to $B(N)$ in $D_N$, while that in Theorem \ref{thm2} may be thought of as  having the information in some sort of blow up of the normal space to $B(N)$ in $D_N$.

    \section{Proof of Theorem III}   \setcounter{equation}{0}
    The proof of Theorem III will be given following several preliminary discussions on the following topics:
 {\setbox0\hbox{(1)}\leftmargini=\wd0 \advance\leftmargini\labelsep
   \begin{itemize}
    \item nilpotent orbits and the reduced limit period mapping;
    \item monodromy cone structure associated to a normal crossing variety;
    \item the differential of the reduced limit period mapping.
    \end{itemize}
    \subsection*{Nilpotent orbits and the reduced limit period mapping}  We begin by recalling some definitions and results from  \cite{CKS1}, \cite{KP1},  \cite{KP2}, \cite{GGK}, \cite{GG} and \cite{GGR}, the last two of whose notations we shall generally follow.  We let
    \begin{itemize}
    \item $\wt B(N)=$ set of nilpotent orbits $(\fbul,N)$.\end{itemize}
    Here, $D=G_\R/H$ is a \mtd\ embedded as an open $G_\R$-orbit in its compact dual $\CD=G_\C/P$.  The \mtd\ structure on $D$ gives a realization of $\CD$ as a set of filtrations $\fbul=\{F^m\subset F^{m-1}\subset\dots\subset F^0=V_\C\}$ on the complexification of a $\Q$-vector space $V$.  The monodromy logarithm  $N\in \cG^\nilp\subset \End(V)$ is a nilpotent endomorphism of $V$ that gives rise to the  monodromy weight filtration, which we center at zero,
    \[
    W_{-k}(N)\subset\dots\subset W_0(N)\subset\cdots\subset W_k(N)=V,\qquad k\leqq m\]
    where $N^{k+1}=0$, $N^k\ne 0$ $(k\leqq m)$.  The conditions to be a nilpotent orbit are
    \begin{enumerate}
    \item $NF^p \subset F^{p-1}$;
    \item $\exp(zN)\cdot \fbul\in D$ for $\rim z\gg 0$.\end{enumerate}}
    \noindent 
    It is  known and of central importance that (\cite{CKS1})
    \begin{multline}\lab{5.1}
    (\fbul,N)\hensp{\em is a nilpotent orbit}\\\iff (V,\wbul(N),\fbul)\hensp{\em is a \lmhs}.\end{multline}
    Here we recall that a \lmhs\ $(V,\wbul(N),\fbul)$ is given by $\fbul$ and $N$   where $\wbul(N)$ is the monodromy weight filtration and where $\fbul$ reduces on $\Gr_k^{\wbul(N)}$ a pure Hodge structure  of weight $k$.  All of our \lmhs s will be \emph{polarized} by a $Q:V\otimes V\to \Q$ (cf.\ \cite{Sc} and \cite{CKS1}).

    Two nilpotent orbits $(N,\fbul)$ and $(N,\fbulp)$ are said to be  \emph{equivalent} if
    \[
    \fbulp = \exp(z N)\cdot \fbul\]
    for some $z\in \C$; i.e., if they lie in the same $\exp(\C N)$ orbit in $\CD$.   We let
    \begin{itemize}
    \item $B(N)=\exp(\C N)\bsl \wt B(N)=$ set of nilpotent orbits modulo equivalence.\end{itemize}
    Assuming that $N\neq 0$ there is a  \emph{reduced \lpm}\ (called a na\"{\i}ve limit in \cite{KP1})
    \begin{equation}\lab{5.2}
    \Phi_\infty:B(N)\to\part D,\end{equation}
    whose image lies in a $G_\R$-orbit.       The definition is
    \[
    \Phi_\infty(\fbul,N)=\lim_{z\to\infty} \exp(zN)\cdot \fbul =: \fbuly.\]
 If we think of $\CD$ as embedded in a product of projective spaces via the Pl\"ucker embeddings of the individual  subspaces $F^p\subset V_\C$, then since $N$ is nilpotent the Pl\"ucker coordinates of $\exp(zN)\cdot F^p$ are polynomials in $z$ and thus have a well-defined limit at $z=\infty$.  In effect  $\Phi_\infty(\fbul,N)$ picks out the highest powers of $z$ in the Pl\"ucker coordinates of $\exp(z N) F^p$.  An elementary general fact is that for any nilpotent $N$ the vector field on $\CD$ induced by the action of the 1-parameter group $\exp(zN)$ vanishes to $2^{\rm nd}$ order at the limit point $\fbuly$, so that the reduced \lpm\ is well defined on the quotient space $B(N)$ of $\wt B(N)$.
   
   One of the important features of the reduced \lpm\ is
   \begin{demo}\lab{new5.3}
   \em The mapping \pr{5.2} factors
   \[
   \xymatrix@R=1.25pc{B(N)\ar[dd] \ar[dr]&\\
   &\part D\\
   B(N)_\R \ar@/_1pc/[uu]\ar@{-->}[ur]&}\]
   through the set $B(N)_\R$ of equivalence classes of $\R$-split \lmhs s.
   \end{demo} 
{\setbox0\hbox{1}\leftmargini=\wd0 \advance\leftmargini\labelsep
   \begin{itemize}
   \item Associated to a \mhs\ $(V,\wbul,\fbul)$ there is the canonical \emph{Deligne bigrading}
   \begin{equation}\lab{new5.4}
   V_\C=\oplus I^{p,q}\end{equation}
   where
   \[
   \bcs
   F^p= \opplus_{r\geqq p}I^{r,\bullet}\\ \\[-11pt]
   W_k=\opplus_{p+q\leqq k} I^{p,q}\\ \\[-11pt]   I^{p,q}\equiv \ol I^{q,p}\hensp{modulo}W_{p+q-2};\ecs\]
   \item The \mhs\ is \emph{$\R$-split} in case
   \[
   I^{p,q}=\ol{I^{q,p}};\]
   canonically associated to a \mhs\ $(V,\wbul,\fbul)$ is an $\R$-split \mhs\  $(V,\wbul,\tfbul)$;
   \item If $(V,\wbul(N),\fbul)$ is a \lmhs, then so is the $\R$-split \mhs\ $(V,\wbul(N),\tfbul)$, and conversely.
     \end{itemize}
     \noindent It follows from this last propery that we have the factorization \pr{new5.3}.  For the time being we will assume that
     \begin{equation}\lab{new5.5}
     (V,\wbul(N),\fbul)\hensp{\em is $\R$-split.}\end{equation}
     In this case the filtration $\fbuly$ is related to $\fbul$ by
     \[
     F^p_\infty = \opplus_{q\leqq m-p} I^{\bullet,q},\]
     where $m$ is the weight of the Hodge structure under consideration.
     \begin{itemize}
     \item If $(V,\wbul,\fbul)$ is a \mhs, then the inclusion $\cG\subset \End_Q(V,V)$ induces on $\cG$ a \mhs\ $(\cG,W_{\bullet,\cG},\fbul_{\cG})$;
     \item Under the assumption \pr{new5.5} we have
     \[
     \cGC = \oplus I^{p,q}_\cG\]
     and
     \begin{equation}\lab{new5.6}
     F^p_{\cG,\infty} = \opplus_{q\leqq p} I^{\bullet,q}_\cG;\end{equation}
     \item the monodromy logarithm $N\in I^{-1,-1}_\cG$;
     
     \item with the identifications
     \[
     \bcs
     T_{\fbul_\cG}\CD= \opplus_{p\leqq -1} I^{p,\bullet}_\cG=F^{-1}_\cG\\ \\[-11pt]
     T_{F^\bullet_{\cG,\infty}} \CD= \opplus_{q\geqq  1} I^{\bullet,q}_\cG =F^{-1}_{\cG,\infty};\ecs\]
     
     \end{itemize}}
     
     \begin{demo}\lab{new5.7}
     \em   the differential
     \[
     \Phi_{\infty,\ast} : T_{\fbul} B(N)_\R\to T_{\fbuly}\CD\]
     of the reduced \lpm\ is the identity on $I^{p,q}_\cG$ for $q\geqq 1$ and is zero on $I^{p,q}_\cG$ for $q< 0$.\end{demo}
     \noindent Pictorially, we picture $I^{p,q}_\cG$ in the $(p,q)$ plane
     \[
     \begin{picture}(100,100)
     \put(0,50){\line(1,0){100}}
     \put(50,100){\line(0,-1){100}}
     \put(105,50){$p$}
     \put(10,80){II}
     \put(85,75){I}
     \put(60,90){$q$}
     \put(10,29){III}\end{picture}\]
 Then
     \begin{align*}
     T_{\fbul}\wt B(N)&\subseteq \hbox{I }\cup \hbox{ II}\\
     T_{\fbuly} \CD&\cong \hbox{II  }\cup \hbox{ III}\end{align*}
     and $\Phi_{\infty,\ast}$ is the identity on the interior of  II with
     \[
     \bcs \phantom{\rm co}\!\ker \Phi_{\infty,\ast}= \hbox{ I}  \\
  \hbox{coker} \,\Phi_{\infty,\ast}=\hbox{ III}.\ecs\]
     \subsection*{Monodromy cone structure associated to a normal crossing variety}
   More generally, associated to a nilpotent cone
   \[
   \sigma=\rspan_{\Q\geqq 0}\{N_1,\dots,N_\ell\}\]
with interior $\sigma^\circ$, from \cite{CKS1}, \cite{CKS2} there is an intricate and deep structure of nilpotent orbits, or equivalently \lmhs s in several variables.  Among the properties of this structure are
   \begin{itemize}
   \item the monodromy weight filtration is independent of $N\in\sigma^0$ (\cite{CK});
   \item the \lmhs\ associated to a nilpotent orbit $(\fbul,N)$ is independent of $N\in\sigma^\circ$;
   \item denoting by $\Delta^\ast(r)$ a punctured disc of radius $r$, in the manner described in \cite{CKS1}  on $\Delta^\ast(r_1)\times\cdots \times \Delta^\ast(r_\ell)$ there are several variable nilpotent orbits  
   \[
   \exp(z_1 N_1+\cdots + z_\ell N_\ell)\cdot \fbul, \qquad \rim z_i\gg 0\]
   which induce variations of mixed Hodge structure (\cite{St-Zu}) on the axes in $\Delta(r_1)\times\cdots\times \Delta(r_\ell)$.\end{itemize}

An important example of this cone structure is provided by a normal crossing variety $X$ for which there exists a $\xi\in T_X\Def(X)$ such that
\begin{equation}\lab{5.3}
\xi_{D_a}\ne 0\hensp{\em for each of  the connected components} D_a, a\in A, \hensp{\em of} D.\end{equation}
  We will describe $\sigma_X$  when the following special condition  is satisfied:
   \begin{align}\lab{5.4}
& \hbox{  \em for each $a\!\in\! A$, there is   $\xi_A\!\in\!\EExt^1_{\cOX}(\Om^1_X,\cO_X)$ such that}\\
&\hbox{ $\xi_{A,a}\!\neq\! 0$ \em while $\xi_{A,b}=0$ for $b\not\in A$.}\notag\end{align}
Geometrically, to first order we may deform $X$ smoothing the component $D_a$ of $X_\sing$ while remaining locally equisingular along the other components $D_b$, $b\neq a$.  Under the assumption \pref{5.4} we may to \fst\ independently smooth the components $D_a$ of the singular locus $D$ of $X$.  Then  there are monodromy transformations $N_a$, $a\in N$, that lead to a nilpotent cone.   

In general the map
\begin{equation}\lab{5.5}
\EExt^1_\cOX(\Omox)\to \opplus_{a\in A}H^0(\cO_{D_a})\end{equation}
will fail to be surjective and additional constructions are needed to obtain a set of monodromy cones described by the combinatorics of how the image of the mapping \pr{5.5} meets the ``coordinate axes'' given by the right-hand term. The details of   this   will be given in the aforementioned work in progress.  In that work we hope to also give the description of the cone in case $X$ is locally a product of normal crossing varieties.  In this situation the faces of the cone will correspond to where factors in the local product description become smoothed, as well as to where components in the stratification of $X_\sing$ become smoothed.

\begin{proof}[Proof of Theorem {\rm \ref{thm3}}]
The proof now follows from a very particular case of Robles' result \cite{Ro}.  In this special case the argument is much simpler and goes as follows.

The first step is to identify the tangent spaces to the $G_\R$-orbit
\[
\cO_{\fbuly}=:G_\R\cdot \fbuly \subset \part D.\]
This is done   in \cite{KP1} and later in \cite{GGK} and \cite{GG}; we shall follow the notations and indexing in the latter.  From Section III.A we have for the real tangent space
\begin{align*}
T^\R_{\fbuly}\cO_{\fbuly}&= \mathrm{Image}\lrc{ \cGR\to \cGC/F^0_\infty \cGC}\\
&\cong \bigoplus_{q>0\atop p\leqq 0} \lrp{\cG^{p,q}+q^{q,p}} \cap \cGR \oplus \bigoplus_{q\geqq p>0}  \lrp{\cG^{p,q}\oplus q^{q,p}} \cap \cGR\\
&\cong \Res_{\C/\R} \Big\{\bigoplus_{q>0\atop p\leqq 0} \cG^{p,q}\Big\}\oplus \bigoplus_{q\geqq p>0}\lrp{\cG^{p,q}\oplus q^{q,p}} \cap \cGR\end{align*}
where $\Res_{\C/\R}$ is the restriction of scalars from $\C$ to $\R$ that maps a complex vector space to the same space now considered as a vector space over $\C\subset \R$.  From the discussion above we see that $T_{(X,\xi)}\Def(X,\xi)$ maps to the first factor, which is in the tangent space to the $G_\R$-orbit $\cO_{\fbuly}$.
\end{proof}

As noted above, what one would like is to show that the interiors of the faces of the cone also map under the reduced \lpm\ to $G_\R$-orbits that are in the closure of the image of $\sig_X$.

      \section{The hierarchy of mixed Hodge structures}   
   \setcounter{equation}{0}
   In this discussion we will  restrict to a standard family  $\cX\to \Delta$.  To this situation there are naturally associated four mixed Hodge structures:
   \begin{enumerate}
   \item the \mhs\ on $ H^m(X)$;
   \item that part of the \lmhs\ that can be constructed from $X$ alone;
   \item the \lmhs\ associated to $\cX\to\Delta$, modulo the equivalence $\fbul\sim \exp(zN)\cdot \fbul$ arising from a change of parameter on $\Delta$;\footnote{We have seen that this data will depend only on the \fst\ neighborhood of $X$ in $\cX$.} \item the \lmhs\ associated to the pair $(X,\xi)$ where $\xi\in T_X\Def(X)$.  
      \end{enumerate}   
      
      We recall from the introduction the    \setcounter{thm}{3} 
   \begin{thm} 
   In a manner to be explained in the proof, there are strict implications
   \[
   {\rm (iv)}\implies {\rm (iii)}\implies {\rm (ii)}\implies {\rm (i)}.\]
   \end{thm} \setcounter{thm}{5}
   Intuitively there is successively  strictly less information in the data arising from the situations (iv), (iii), (ii), (i).  Before turning to the proof we will illustrate the result in the following  
   
   \begin{exam}\hspace*{-5pt}\footnote{This is an extension to $g=3$ of the case $g=2$ in \cite{Ca}.}
   Suppose that $X_t$ is a smooth curve of genus $g=3$ whose limit $X$ is an irreducible nodal  curve whose normalization $\wt X$ has genus $\wt g=1$
   \[
 \begin{picture}(250,145)
\put(0,29){\includegraphics[scale=.8]{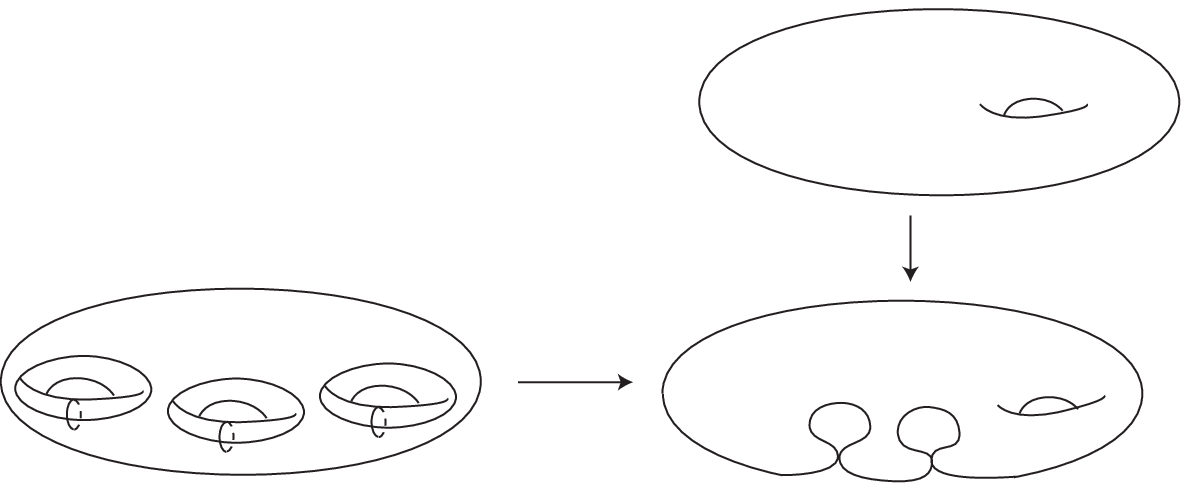} }      
\put(15,25){$\delta_1$}
\put(55,25){$\delta_2$}
\put(85,25){$\delta_3$}\put(25,67){$\gamma_1$}
\put(57,67){$\gamma_2$}
\put(92,64){$\gamma_3$}
\put(187,57){$\gamma_1$}
\put(207,57){$\gamma_2$}
\put(280,50){$X$}
\put(280,120){$\widetilde X$}
\put(178,125){${\scriptstyle \bullet}\ p_1$}
\put(178,110){${\scriptstyle \bullet}\ q_1$}
\put(204,125){${\scriptstyle \bullet}\ p_2$}
\put(204,110){${\scriptstyle \bullet}\ q_2$}
\put(-30,50){$X_t$}   \end{picture}\vspace*{-22pt}\]  
Setting $\ell(t) = (1/2\pi i)\log t$ and using the   symplectic basis drawn above for $H_1(X_t,\Z)$, the normalized period matrix is   
\[
\begin{pmatrix}
1&\\
&1\\
&&1\\
\ell(t)+a_{11}(t)&a_{12}(t)&b_1 (t)\\
a_{21}(t)&\ell(t)+a_{22}(t)& b_2(t)\\
b_1(t)&b_2(t)&c(t)
\end{pmatrix}\]
where the $a_{ij}(t)=a_{ji}(t)$, $b_i(t)$ and $c(t)$ are holomorphic in the disc and $\rim c(t)>0$.  With the choice $t$ of parameter the nilpotent orbit is
\[
\bpm
1\\
&1\\
&&1\\
\ell(t)+a_{11} &a_{12}&b_2\\
a_{21}&\ell(t)+a_{22}&b_2\\
b_1&b_2&c\epm\]
where $a_{ij}=a_{ij}(0)$, $b_i=b_i(0)$ and $c=c(0)$.  Letting $\om_i(t)$ be the holomorphic differentials on $X$ with limits $\om_i$ on $X$ that pull up to $\wt \om_i$ on $\wt X$, we have
\begin{itemize}
\item $\wt \om_1,\wt \om_2$ are differentials of the $3^{\rm rd}$ kind on $\wt X$ with divisor $p_i+	q_i$ and $\Res_{p_i}\wt \om_i=+1$, $\Res_{q_i}\wt \om_i=-1$ for $i=1,2$;
\item $\wt \om_3$ is a holomorphic differential on $\wt X$.\end{itemize}
Under a reparamatrization $t'=e^{2\pi i\la}t$,
\[
a_{ii}(t')=a_{ii}(t)+\la,\]
and all other entries in the period matrix evaluated at $t=0$ are unchanged.

We note that implicit in the choice of symplectic basis is the monodromy weight filtration
\[
\bcs
\hfill W_0&\chs=\rspan\{\delta_1,\delta_2\}\\
W_2/W_1&\chs\cong \rspan\{\gamma_1,\gamma_2\}\\
W_1/W_0&\chs\cong \rspan \{\delta_3,\gamma_3\}.\ecs\]

The entries in the above  period matrix at $t=0$ are (cf.\ \cite{Ca}) \begin{itemize}
\item[($I_1$)] $c$ is the period of the elliptic curve $\wt X$;
\item[($I_2$)] $b_i$ is the image of $\AJ_{\wt X}(p_i-q_i)$ in $J(\wt X)$;
this gives the extension data in
\[
0\to \Gr_0\to \Gr_1\to \Gr_1/\Gr_0\to 0\]
as described in \cite{Ca};
\item[($I_3$)] with suitable normalization of the $\wt \om_i$,
\[
a_{ij}=\int^{p_j}_{q_i}\wt \om_i,\qquad i\neq  j;\]
\item[($I_4$)] finally,  with a choice of parameter $t$ we may uniquely define the improper integrals
\[
\int^{p_i}_{q_i}\wt\om_i.\]
\end{itemize}
This means that if $\gamma_{i,t}$ is the above curve on $X_t$ for $t\neq  0$
\[
\int_{\gamma_{i,t}}  \om_{i,t}=\ell(t) +a_{ii}(t),\]
and then on $\wt X$ we will have
\[
\lim_{\left\{\begin{smallmatrix} p'_i\to p_i\\ q'_i\to q_i\end{smallmatrix}\right.} \int^{p'_i}_{q'_i}\wt \om_i=a_{ii}(0)\]
 where the picture is
\[
\begin{picture}(100,50)
\put(0,25){\includegraphics[scale=.5]{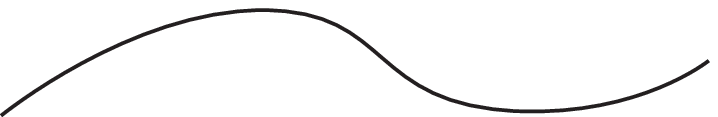}}
\put(-1,24){$\scriptstyle\bullet$}
\put(15,33){$\scriptstyle\bullet$}
\put(15,25){$\scriptstyle p'_i$}
\put(77,24){$\scriptstyle \bullet$}
\put(77,14){$\scriptstyle q'_i$}
\put(100,30){$\scriptstyle \bullet$}
\put(98,24){$\scriptstyle q_i$}
\put(-3,18){$\scriptstyle p_1$}\end{picture}\]
and where the logarithmic singularities at the endpoints cancel.

The $a_{ij}$ for $i\neq j$ record the part of the ``extension upon extension" data in $\Gr_2/G_1$ over $\Gr_1/\Gr_0$ that is invariant under reparametrization, and the $a_{ii}$ record the full extension data.

Algebro-geometrically the picture is the following.  Denoting by $\cM_g$ the moduli space of the stable curves of genus $g$, and by $\ol\cM_g$ the Deligne-Mumford compactification, the curve $X$ gives a point in 
\[
\part {\cM}_3\subset \ol{\cM}_3.\]
More specifically, $X$ defines a point in a codimension-1 component $\cC$ of the stratified variety $\part{\cM}_3$.  Then
\begin{itemize}
\item $\dim \cC=4$ and $c,b_1,b_2$ and $a_{12}$ are local coordinates in $\cC$;
\item $a_{11}$ and $a_{22}$ give normal parameters to $\cC$ in $\ol \cM_3$.\end{itemize}

The difference between (ii) and (iii) in this case is that (iii) contains the information in the weight filtration, which is information that is not obtainable from that on $X$ alone (see the subsequent discussion).  \end{exam}

To explain (ii) we will picture a \lms\ in terms of the $N$-strings as
\[
\begin{matrix}
\qquad H^0(-m)\to\enspace \cdot\enspace \cdot\enspace\cdot\enspace\cdot\enspace \to H^0(-1)\xri{N}H^0\\
 \vdots \\
 H^{m-1}(-1)\to H^{m-1} \\
 H^m \end{matrix}\]
 where $H^k$ is a pure Hodge structure of weight $k$.  It is this presentation that is especially useful in the computation of examples \cite{GG} and \cite{GGR}.  Equivalent data to the above are the following parts of a polarized \lmhs\ $(V,\wbul(N),\fbul)$
 \begin{itemize}
 \item the Hodge structures $\Gr^{\wbul(N)}_k$;
 \item the iterated $N$ operators on the $N$-strings.  \end{itemize} It is known \cite{GG}, \cite{Ro}   that this data always arises from a non-unique \lmhs.

 Our main result, Theorem \ref{thm6} in the introduction, is that, under the assumption that $X$ is smoothable but with no $\xi\in T^0_X\Def(X)$ singled out, we may compute the $H^{m-i}(-j)$'s above purely in terms of $X$ alone.  For this we will use the maps
 \[
 \hbox{Rest}: H^q \bp{X^{[k]}} \to H^q \bp{X^{[k+1]}}\]
 obtained by the alternating sums of the restriction maps, and the suitably alternated  Gysin maps
 \[
 \Gy: H^q\bp{X^{[k]}} (-1) \to H^{q+2}\bp{X^{[k-1]}}.\]
The result is

 \begin{thm} \lab{thm7} Assuming that $X$ is smoothable, there are complexes
 \[\lower.67in\hbox{$ \begin{matrix}
 \cdot\\[-4pt] \cdot\\[-4pt] \cdot\\[-4pt] \cdot\\[-4pt] \cdot \end{matrix}$}
 \xymatrix@R=.25pc{ &H^{q-4}\bp{X^{[k+2]}}(-2)&&\\
 &\oplus&H^{q-2}\bp{X^{[k+1]}}(-1)&\\
  \ar[r]&H^{q-2}\bp{X^{[k]}}(-1)\ar[r]&\oplus\ar[r]&H^q\bp{X^{[k]}}\\\
  &\oplus&H^q\bp{X^{[k-1]}}&\\
  &H^q\bp{X^{[k-2]}}&}\]
  such that for $0\leqq j\leqq m-i$
  \begin{equation}\lab{6.1}
  H^{m-i}(-j)\cong H^\ast_{\rm Rest} H^\ast_{\Gy}\lrp{H^{m-i}\bp{X^{[i+1]}}}(-j).\end{equation}
 \end{thm}
The notation in \pr{6.1} means that the left-hand side is computed by the cohomology of complexes in the statement of the theorem at the spot designated by the right-hand side in \pr{6.1}.  
  A significant fact here is that 
 \begin{equation}\lab{6.2}
 X\hbox{ \em smoothable } \implies \Rest\circ\Gy=- \Gy\circ \Rest.\end{equation}
 The point is that this 
implication \pr{6.2} is generally not true unless $X$ is smoothable, although  the result does not depend on any particular smoothing.

A further significant point is that in taking the cohomology of the restriction sequences
\begin{align*}
H^\ast_{\Gy}\lrp{ H^i \lrp{X^{[j-1]}} (-m)}&\to H^\ast_{\Gy} \lrp{H^i \bp{X^{[j]}}(-m)}\\& \to H^\ast_{\Gy}
\lrp{H^i \bp{X^{[j+1]}}(-m)}
\end{align*}
we only put in $H^i\bp{X^{[j]}}(-k)$ if $0\leqq k<j-1$.

Referring to \pr{1.9} and Theorem \ref{thm6}   in the introduction, Theorem \ref{thm7} implies that result, and the $N$-maps in the $N$-strings are the twisted identity maps on the individual pieces as given by \pr{6.1}.
  
 For the proof of  Theorem \ref{thm7} one uses the basic constriction introduced in \cite{St1} and \cite{Zu}.  We will recall this for a standard family $\cX\xri{\pi}\Delta$ and observe at the end that the vector space $\bH^m\lrp{\Obul_{\cX/\Delta}(\log X)\otimes \cOX}$ has a filtration shifted down two steps by $N$, and the associated graded together with the mappings induced by $N$ may be defined in terms of $X$ alone.  The $N$-strings that result are the ones expressed in the  theorem.  The upshot is that given an abstract normal crossing variety $X$ we will be able to define the object that is defined in terms of $X$ alone, provided only we assume that
$ X$ is smoothable. 
 This object is related to the \lmhs\ associated to the standard family in the manner just described.  It is interesting to note that the object so described will be \emph{independent} of the smoothing of $X$, provided   that one exists.  This is a reflection of the result in \cite{CKS1} that in a multi-parameter family the \lmhs\ associated to the origin is independent of the direction of approach from  the interior of the cone.
 
 The construction in \cite{St1} and \cite{Zu} goes as follows: We define a bi-graded complex $A^{\bbul}$ where
 \[
 A^{p,q}=\Om^{p+q+1}_\cX(\log X) \big\slash \wt W_q\Om^{p+q+1}_\cX(\log X),\]
 where $\wt W_q$ is the standard filtration given by
 \[
 \wt W_q = \hbox{ differential forms with at most   $q$ $d x_i/x_i$ terms.}\]
 The differentials are given by
 \[
 \bcs d'&\chs= \hbox{usual } d\\
 d''&\chs= \wedge dt/t.\ecs\]
 The basic observation and definitions  are
  {\setbox0\hbox{(1)}\leftmargini=\wd0 \advance\leftmargini\labelsep
\begin{itemize}
 \item The mapping 
 $\Om^p_\cX(\log X)\xri{\wedge dt/t} A^{p,0}$ has co-kernel naturally isomorphic to $\Om^p_{\cX/\Delta}(\log X)\otimes \cOX$; 
   \item $W_k A^{p,q}=:\wt W_{2q+k+1}\Om^{p+q+1}_\cX(\log X) \big\slash \wt W_q \Om^{p+q+1}_\cX(\log  X) $;
 \item $F^r A^{p,q} = \bcs A^{p,q}&\hbox{if }p\geqq r,\\
 0&\hbox{if } p<r.\ecs$ 
 \end{itemize}}%
 \noindent 
 Then with respect to $\fbul$ the map
 \[
 \Obul_{\cX/\Delta}(\log X) \otimes \cOX\to \abul\]
 is a filtered quasi-isomorphism. By definition $\wbul\abul$ induces the weight filtration on $\bH^m(\Obul_{\cX/\Delta}(\log X))\otimes \cOX$.  We note that
 \begin{align*}
 \Gr^W_k A^{p,q}& = W_{k+2q+1}\Om^{p+q+1}_\cX(\log X) \big\slash W_{k+2q}\Om^{p+q+1}_\cX(\log X)\\
 &\cong \bcs 0&\hbox{if }k\leqq -(q+1),\\
 \Om^{p-q-k}_{X^{[2q+k+1]}}&\hbox{if }k\geqq -q\ecs\end{align*}
 where the second isomorphism is by the iterated residue map.
 
 An intermediate step to computing $\bH^\ast(A)$ is to use the spectral sequence associated to $\wt \wbul$.  For this spectral sequence where $A^i=\opplus_{p+q=i}A^{p,q}$\pagebreak
 \begin{align} \lab{6.6}
 E^{a,b}_1& = H^{a+b}\lrp{\Gr^W_{-b}A^i}\\
 &=H^{a+b} \lrp{\opplus_q \Om^{i-2q+b}_{X^{[2q+1-b]}}}\notag\\\notag
 &=\bigoplus_{q=\min (0,b)} H^{i-2q+b,a+b}\lrp{X^{[2q+1-b]}}.\end{align}
 If $\dim X=n$, then $\dim X^{[2q+1-b]}= n+1-(2q+1-b) = n+b-2q$, from which we have
 \[\displaylines{
 H^{i-2q+b,a+b}\lrp{X^{[2q+1-b]}}\neq 0\hfill\cr\hfill\implies \bcs
 i-2q+b\leqq n-2q+b=\dim X^{[2q+1-b]}\iff i\leqq n\\
 a+b\leqq n-2q+b\iff a\leqq n-2q.\ecs}\]
 These are the only inequalities other than $\max(a,b)\leqq q\leqq i$.  Thus to have potentially non-zero $H^{r,s}\bp{X^{[t]}}$ we need 
 \[
 \begin{array}{l}
 i-2q+b=r\\
 a+b=s\\
 t=2q+1-b\end{array}\iff \begin{array}{l}
 i=t+r+1\\
 b=2q+1-t\\
 a=s+t-2q-1\end{array} (\implies b\equiv 1-t(\mod 2))\]
 and $q\geqq \max(0,b)$ gives $b\geqq \max(1-t,2b+1-t)$, which then  gives
 \[
 \bcs
 t-1&\chs \geqq b\geqq 1-t\\
 \hfill b&\chs\equiv 1-t(\mod 2)\ecs\]
 (thus $t=1,b=0$; $t=2,b=-1,1$; $t=3, b=-2,0,2,\ldots$) for $r=i-t+1=i+\dim X^{[t]}-n=\dim X^{[t]}-(n-i)$.  Fixing $t$ and noting $0\leqq i\leqq n$ we then have
 \[
 H^{r,\ast}\bp{X^{[t]}}\hbox{ appears for one value of $i$ and $t$ values of $b$.}\]
 
  The above gives the conclusion that for ${H^n}={\bH^n}(\hbox{LMHS})$, before cancellation in the spectral sequence 
 \begin{alignat*}{5}
 {H^n}\bp{X^{[1]}}& \hensp{appears once}&& (b=0)\\
  { H^{n-1}}\bp{X^{[2]}}&\hensp{appears once} &&(b=-1,1)\\
 { H^{n-2}}\bp{X^{[3]}}&\hensp{appears three times}&& (b=-2,0,2)\\
 &\vdots\end{alignat*}
 
 \emph{Sketch of the proof of Theorem {\rm \ref{th4}}.}
 Referring now to \pr{6.6}, by \cite{St1} the spectral sequence degenerates at $E_2$.  The $E_1$-terms are as indicated there.  The $d_1$-map is, as noted in \cite{Zu}, ``composed of various restriction  and Gysin maps.''  After unwinding the indices, the $d_1$-complex turns out to be the one in the statement of the theorem.
 
 We now come to the main point.  Note that the individual terms and maps may be defined in terms of $X$ alone.  However in general we do \emph{not} have
 \begin{equation}\lab{6.8}
 \Rest\circ \Gy = - \Gy\circ \Rest.\end{equation}
 
 \begin{propn}
 The anti-commutativity  commutativity \pr{6.8} holds if, and only if, $\cO_D(X)$ is topologically trivial.  This is the case if $X$ is smoothable.
 \end{propn}
 
 Rather than give the formal argument we shall indicate by example in the simplest non-trivial cases why the result should be true.\footnote{The details of this argument will appear in the previously mentioned work in progress.}  We note that if $X$ is smoothable  the proposition is true.  Our central point is that this  sufficient condition is essentially also necessary.  We say essentially, because \pr{6.8} is a purely topological fact which only requires that $\cO_D(X)$ be topologically, but not necessarily analytically, trivial.\looseness=-1
 
 Let $X$ be an irreducible surface having as singular locus a double curve $C$ whose inverse image in $X^{[1]}$  is a disjoint union $X^{[1]}= C_1\amalg C_2$ of two smooth curves.  We will denote by $H^q(X^{[2]})_-$ the classes $\alpha\oplus  -\alpha \in H^q(C_1\amalg C_2)$.  Then $H^q(X^{[2]})_-\cong H^q(C)$, but we put opposite signs on those in $H^q(C_1)$ and $H^q(C_2)$.  Then we shall show \vspth
 \begin{enumerate}
 \item $\cO_D(X)$ is topologically trivial if, and only if, $C^2_1=-C^2_2$;
 \item the sequence
 \begin{equation}\lab{6.9}
 H^0\bp{X^{[2]}}_- (-1)\xri{\Gy} H^2\bp{X^{[1]}}\xri{\Rest} H^2\bp{X^{[2]}}_-\end{equation}
 is a complex if, and only if, $C^2_1=-C^2_2$.\end{enumerate}
 The complex \pr{6.9} is the simplest non-trivial case of the sequences that appear in the statement of Theorem \ref{thm7}.  
 
 Denoting by $\eta_{C_i}\in H^2\bp{X^{[1]}}$ the fundamental class of $C_i$ and by $[C_i]$ the fundamental class of $C_i$ itself, the sequence \pr{6.9} is
 \[
 \begin{array}{lc}
 1_{C_1}-1_{C_2} \to \eta_{C_1}-\eta_{C_2}\longrightarrow& \lrc{\bp{(C_1-C_2)\cdot C_1}[C_1]+\bp{(C_1-C_2)\cdot C_2[C_2]}}\\
 &\sideeq\\
& \lrc{C^2_1[C_1]-C^2_2[C_2]}_-\\
& \sideeq\\
& \lrp{\frac{C^2_1+C^2_2}{2}} \bp{[C_1]-[C_2]},\end{array}\] 
 which proves (ii).
 
 As for (i) we have
\vspth \[
 N_D(X)\cong N_{C_1}\bp{X^{[1]}} \otimes N_{C_2}\bp{X^{[1]}}^\ast.\vspth\]
 
 In case $X$ is still an irreducible surface a piece of the complex in Theorem \ref{thm7} is
 \vspth\[
 \xymatrix{
 &H^0\bp{X^{[3]}}_-\ar[dr]^u&\\
 H^0\bp{X^{[2]}}_-(-1)\ar[dr]_g\ar[ur]^f &\oplus& H^2\bp{X^{[2]}}_-\\
 &H^2\bp{X^{[1]}}\ar[ur]_v&}\vspth\]
 Here the minus sign on $\bp{\quad}_-$ refers to classes that transform by the sign of the induced action on cohomology given the labeling into even-odd of the components lying over a general point in the map $X^{[k]}\to X_k$.  The conditions to have a complex are
 \[
 u\cdot f+v\cdot g=0,\]
 which when worked out is a consequence of $\cO_D(X)\cong \cO_D$.

\end{document}